\documentclass{article}

\usepackage{proof}
\usepackage{latexsym}
\usepackage{amsfonts}
\usepackage{amssymb}

\usepackage{color}

\newcommand{\blem}{\begin{lemma}}
\newcommand{\elem}{\end{lemma}}
\newcommand{\bth}{\begin{theorem}}
\newcommand{\ethm}{\end{theorem}}
\newcommand{\benu}{\begin{enumerate}}
\newcommand{\eenu}{\end{enumerate}}
\newcommand{\bdes}{\begin{description}}
\newcommand{\edes}{\end{description}}
\newcommand{\bdf}{\begin{definition}}
\newcommand{\edf}{\end{definition}}
\newcommand{\bcor}{\begin{cor}}
\newcommand{\ecor}{\end{cor}}
\newcommand{\bprp}{\begin{proposition}}
\newcommand{\eprp}{\end{proposition}}
\newcommand{\bmlem}{\begin{mlemma}}
\newcommand{\emlem}{\end{mlemma}}
\newcommand{\bclm}{\begin{claim}}
\newcommand{\eclm}{\end{claim}}
\newcommand{\bprf}{{\bf Proof}.\hspace{2mm}}

\newcommand{\eprf}{\hspace*{\fill} $\Box$}

\newcommand{\beqn}{\begin{equation}}
\newcommand{\eeqn}{\end{equation}}
\newcommand{\beqnarr}{\begin{eqnarray}}
\newcommand{\eeqnarr}{\end{eqnarray}}
\newcommand{\beqnarrs}{\begin{eqnarray*}}
\newcommand{\eeqnarrs}{\end{eqnarray*}}

\newtheorem{theorem}{Theorem}[section]
\newtheorem{definition}[theorem]{Definition}
\newtheorem{proposition}[theorem]{Proposition}
\newtheorem{lemma}[theorem]{Lemma}
\newtheorem{cor}[theorem]{Corollary}
\newtheorem{mlemma}[theorem]{Main Lemma}
\newtheorem{claim}[theorem]{Claim}

\title{
Proof theory of weak compactness
}
\author{Toshiyasu Arai
\\
Graduate School of Science,
Chiba University
\\
1-33, Yayoi-cho, Inage-ku,
Chiba, 263-8522, JAPAN
\\
tosarai@faculty.chiba-u.jp
}
\date{}
\begin{document}
\maketitle

\begin{abstract}
We show that the existence of a weakly compact cardinal over the Zermelo-Fraenkel's set theory ${\sf ZF}$
is proof-theoretically reducible to iterations of 
Mostowski collapsings and 
Mahlo operations.
\end{abstract}


\section{Introduction}\label{sect:introduction}

It is well known that a cardinal is weakly compact iff it is $\Pi^{1}_{1}$-indescribable.
From this characterization we see readily that the set of Mahlo cardinals below a weakly compact cardinal
is stationary, i.e., every club (closed and unbounded) subset of a weakly compact cardinal contains a Mahlo cardinal.
In other word, any weakly compact cardinal is hyper Mahlo.
Furthermore any weakly compact cardinal $\kappa$ is in the diagonal intersection $\kappa\in M^{\triangle}=\bigcap\{M(M^{\alpha}):\alpha<\kappa\}$
for the $\alpha$-th iterate $M^{\alpha}$ of the Mahlo operation $M$:
for classes $X$ of ordinals,
\[
\kappa\in M(X) :\Leftrightarrow
X\cap\kappa \mbox{ is stationary in } \kappa
\Leftrightarrow
\forall Y\subset \kappa[(Y \mbox{  is club}) \to X\cap Y\neq\emptyset].
\]
Note that $\kappa\in M(X)$ is $\Pi^{1}_{1}$ on $V_{\kappa}$.

On the other side R. Jensen\cite{Jensen} showed under the axiom $V=L$ of constructibility 
that for  regular cardinals $\kappa$,
$\kappa$ is weakly compact  iff
$\forall X\subset\kappa[\kappa\in M(X)\Rightarrow M(X)\cap\kappa\neq\emptyset]$
iff
$\forall X\subset\kappa[\kappa\in M(X)\Rightarrow \kappa\in M(M(X))]$.

Jensen's proof in \cite{Jensen} yields a normal form theorem of $\Pi^{1}_{1}$-formulae on $L_{\kappa}=J_{\kappa}$
uniformly for regular uncountable
 cardinals $\kappa$ as follows.

For a first order formula $\varphi[D]$ with unary predicates $A,D$,
let 
 \begin{eqnarray}
&& \alpha\in S^{\varphi}(A)
 :\Leftrightarrow \mbox{ there exists a limit } \beta \mbox{ such that } \alpha<\beta<\alpha^{+}, 
 A\cap\alpha\in J_{\beta},
 \nonumber
 \\
 &&  \langle J_{\beta},\in,A\cap\alpha\rangle\models 
 \forall D\subset\alpha\, \varphi[D], \alpha \mbox{ is regular in } \beta
 \mbox{ and } 
 \nonumber
 \\
 &&
 \exists p\in J_{\beta}\forall X[(p\cup\{\alpha\}\subset X\prec J_{\beta}) \land (X\cap\alpha\mbox{ is transitive}) \Rightarrow X=J_{\beta}]
 \label{eq:Jensenmin}
 \end{eqnarray}
where $\alpha$ is regular in $\beta$ iff there is no cofinal function from a smaller ordinal$<\alpha$
into $\alpha$, which is definable on $J_{\beta}$.

The following Proposition \ref{prp:JensenLemma5.2} is the Lemma 5.2 in \cite{Jensen}.

\begin{proposition}\label{prp:JensenLemma5.2}
Let $\alpha\in S^{\varphi}(A)$ and $\beta$ be an ordinal as in the definition of $S^{\varphi}(A)$.
Then $\alpha$ is $\Sigma_{1}$-singular in $\beta+1$, i.e.,
there exists a cofinal function from a smaller ordinal$<\alpha$
into $\alpha$, which is $\Sigma_{1}$-definable on $J_{\beta+1}$.
\end{proposition}

Fix a regular uncountable cardinal $\kappa$, a set $A\subset\kappa$.
For a finite set $\{A,\ldots\}$ of subsets $A,\ldots$ of $\kappa$
and ordinals $\alpha<\kappa$, let $N_{\alpha}(A,\ldots)$ denote the least $\Sigma_{1}$-elementary submodel of $J_{\kappa^{+}}$,
$N_{\alpha}(A,\ldots)\prec_{\Sigma_{1}}J_{\kappa^{+}}$, such that $\alpha\cup \{A,\ldots\}\cup\{\kappa\}\subset N_{\alpha}(A,\ldots)$.
Namely $N_{\alpha}(A,\ldots)$ 
is the $\Sigma_{1}$-Skolem hull $\mbox{Hull}^{J_{\kappa^{+}}}_{\Sigma_{1}}(\alpha\cup \{A,\ldots\}\cup\{\kappa\})$ of 
$\alpha\cup \{A,\ldots\}\cup\{\kappa\}$ on $J_{\kappa^{+}}$.
Let
\[
C(A,\ldots):=\{\alpha<\kappa: N_{\alpha}(A,\ldots)\cap\kappa\subset\alpha\}
.\]
Then it is easy to see that $C(A,\ldots)$ is club in $\kappa$, and definable over $J_{\kappa^{+}}$.

\begin{proposition}\label{prp:Jensen}
Let $\kappa$ be a regular uncountable
 cardinals $\kappa$, $A\subset\kappa$, $\varphi[D]$ a first order formula with parameters $A,D$.
\benu
 \item\label{prp:Jensen1}
Suppose $\langle J_{\kappa^{+}},\in,A\rangle\models  \forall D\subset\kappa\, \varphi[D]$, and let $C$ be a club subset of $\kappa$.  
Then  the least element of the club set $C(A,C)$ is in $S^{\varphi}(A)$.
 \item\label{prp:Jensen2}
Suppose $\langle J_{\kappa^{+}},\in,A\rangle\not\models  \forall D\subset\kappa\, \varphi[D]$, then 
  $S^{\varphi}(A)\cap C(A)=\emptyset$.
\eenu

 Thus
$\langle J_{\kappa^{+}},\in,A\rangle\models  \forall D\subset\kappa\, \varphi[D]$
iff $S^{\varphi}(A)$ is stationary in $\kappa$.
And $\kappa$ is weakly compact iff
for any stationary subset $S\subset\kappa$ there exists an uncountable regular cardinal $\alpha<\kappa$
such that $S\cap\alpha$ is stationary in $\alpha$.
\end{proposition}
\bprf

 \ref{prp:Jensen}.\ref{prp:Jensen1}.
Suppose $\langle J_{\kappa^{+}},\in,A\rangle\models  \forall D\subset\kappa\, \varphi[D]$, and let $C$ be a club subset of $\kappa$.
Consider the 
club subset $C(A,C)$ of $\kappa$.
Then $C(A,C)\subset C$ .
We show that $\alpha\in S^{\varphi}(A)$ for the least
element $\alpha$ of $C(A,C)$.
Let $\pi:\langle J_{\beta},\in , A\cap\alpha,C\cap \alpha\rangle\cong N_{\alpha}(A,C)\prec_{\Sigma_{1}}J_{\kappa^{+}}$ be the transitive collapse of $N_{\alpha}(A,C)$.
$\beta$ is a limit ordinal with $\alpha<\beta<\alpha^{+}$.
From $\langle J_{\kappa^{+}},\in,A\rangle\models  \forall D\subset\kappa\, \varphi[D]$ we see 
$\langle J_{\beta},\in , A\cap\alpha\rangle\models  \forall D\subset\alpha\, \varphi[D]$, and $A\cap\alpha,C\cap\alpha\in J_{\beta}$ from $A,C\in N_{\alpha}(A,C)$.
It remains to show (\ref{eq:Jensenmin}) for $p=\{A\cap\alpha,C\cap\alpha\}$.
Assume $\{A\cap\alpha,C\cap\alpha,\alpha\}\subset X\prec J_{\beta}$ 
and $X\cap\alpha=\gamma$ for an ordinal $\gamma\leq\alpha$.
Then $\gamma\cup\{A,C,\kappa\}\subset\pi"X\prec N_{\alpha}(A,C)\prec_{\Sigma_{1}} J_{\kappa^{+}}$.
This yields $N_{\gamma}(A,C)\prec_{\Sigma_{1}}\pi"X$, and 
$N_{\gamma}(A,C)\cap\kappa\subset(\pi"X)\cap\kappa=\pi"(X\cap\alpha)=\gamma$ by $N_{\alpha}(A,C)\cap\kappa\subset\alpha$.
This means that $\gamma\in C(A,C)$, and hence $X\cap\alpha=\gamma=\alpha$.
Therefore $\pi"X=N_{\alpha}(A,C)$, and 
$X=J_{\beta}$.
\\

\noindent
\ref{prp:Jensen}.\ref{prp:Jensen2}.
Suppose $\langle J_{\kappa^{+}},\in,A\rangle\not\models  \forall D\subset\kappa\, \varphi[D]$.
Assume $\alpha\in S^{\varphi}(A)\cap C(A)$.
Let $\langle J_{\bar{\beta}},\in , A\cap\alpha\rangle\cong N_{\alpha}(A)\prec_{\Sigma_{1}}J_{\kappa^{+}}$ be the transitive collapse of $N_{\alpha}(A)$.
Then $\langle J_{\bar{\beta}},\in , A\cap\alpha\rangle\not\models  \forall D\subset\alpha\, \varphi[D]$.
On the other hand we have by $\alpha\in S^{\varphi}(A)$, there exists a limit $\beta$ such that
$\langle J_{\beta},\in , A\cap\alpha\rangle\models  \forall D\subset\alpha\, \varphi[D]$,
and $\alpha$ is $\Sigma_{1}$-singular in $\beta+1$ by Proposition \ref{prp:JensenLemma5.2}.
Hence $\beta<\bar{\beta}$ and $\alpha$ is $\Sigma_{1}$-singular in $\bar{\beta}$.
This means that $\kappa$ is $\Sigma_{1}$-singular in $\kappa^{+}$.
However $\kappa$ is assumed to be regular.
A contradiction.
\eprf

In this paper we show that 
the existence of a weakly compact cardinal over the Zermelo-Fraenkel's set theory ${\sf ZF}$
is proof-theoretically reducible to iterations of 
Mostowski collapsings and 
Mahlo operations.

Let $\mathbb{K}$ denote the formula stating that `there exists a weakly compact cardinal $\mathcal{K}$'.

For $\Sigma^{1}_{2}$-sentences $\varphi\equiv\exists Y\forall X\,\theta$, let $\varphi^{V_{\mathcal{K}}}$ be
$\exists Y\subset V_{\kappa}\forall X\subset V_{\mathcal{K}}\,\theta^{V_{\mathcal{K}}}$ where $\theta^{a}$ denotes the result of restricting any unbounded quantifiers $\exists x,\forall x$ to $\exists x\in a, \forall x\in a$, resp.

\begin{theorem}\label{th:mainthK}
There are $\Sigma_{n+1}$-formulae $\theta_{n}(x)$ for which the following holds:
\benu

\item\label{th:mainthK1}
For each $n<\omega$,
\[
{\sf ZF}+(V=L)\vdash
\forall\mathcal{K}[(\mathcal{K} \mbox{ {\rm is a weakly compact cardinal}})  \to  \theta_{n}(\mathcal{K})]
\]
and
\[
{\sf ZF}+(V=L)\vdash
\forall\mathcal{K}[\theta_{n+1}(\mathcal{K}) \to \mathcal{K}\in M(\{\pi<\mathcal{K}: \theta_{n}(\pi)\})]
.\]

\item\label{th:mainthK2}
For any $\Sigma^{1}_{2}$-sentences $\varphi$, 
if
\[
{\sf ZF}\vdash\forall\mathcal{K}[(\mathcal{K} \mbox{ {\rm is a weakly compact cardinal}}) \to \varphi^{V_{\mathcal{K}}}]
,
\]
then we can find an $n<\omega$ such that
\[
{\sf ZF}+(V=L)\vdash
\forall\mathcal{K}[\theta_{n}(\mathcal{K}) \to
\varphi^{V_{\mathcal{K}}}]
.\]

\eenu

Hence ${\sf ZF}+(V=L)+(\mathcal{K} \mbox{ {\rm is weakly compact}})$ is $\Sigma^{1}_{2}(\mathcal{K})$-conservative over
${\sf ZF}+(V=L)+\{\theta_{n}(\mathcal{K}) : n<\omega\}$, and
${\sf ZF}+(V=L)+\mathbb{K}$ is conservative over 
${\sf ZF}+(V=L)+\{
\exists\mathcal{K}\, \theta_{n}(\mathcal{K}) : n<\omega\}$,
e.g.,
with respect to first-order sentences $\varphi^{V_{I_{0}}}$ for the least weakly inaccessible cardinal $I_{0}$.

Note that $T_{n}={\sf ZF}+(V=L)+\{\exists\mathcal{K}\, \theta_{n}(\mathcal{K})\}$
is weaker than ${\sf ZF}+\mathbb{K}$, e.g., ${\sf ZF}+\mathbb{K}$ proves the existence of a model
of $T_{n}$ for each $n<\omega$.
\end{theorem}
The $\Sigma_{n+1}$-formulae $\theta_{n}(x)$ are defined by
\[
\theta_{n}(x) :\Leftrightarrow x\in Mh_{n}^{\omega_{n}(I+1)}
.\]
The $\Sigma_{n+1}$-class $Mh_{n}^{\xi}$ for ordinals $\xi$ is defined through iterations of Mostowski collapsings and Mahlo operations,
cf. Definition \ref{df:Cpsiregularsm}.

Let us explain some backgrounds of this paper.
$\Pi_{3}$-reflecting ordinals are known to be recursive analogues to weakly compact cardinals.
Proof theory ({\it ordinal analysis\/}) of $\Pi_{3}$-reflection has been done by M. Rathjen\cite{Rathjen94},
and \cite{odpi3, ptpi3, wfnonmon2, esubpi2}.

As observed in \cite{ptpi3, LMPS},
ordinal analyses of $\Pi_{N+1}$-reflection yield a proof-theoretic reduction of $\Pi_{N+1}$-reflection 
in terms of iterations of $\Pi_{N}$-recursively Mahlo operations.
Specifically we show the following Theorem \ref{th:resolutionpi3} in \cite{consv}.
Let {\sf KP}$\omega$ denote the Kripke-Platek set theory with the axiom of Infinity, 
$\Pi_{N}(a)$ a universal $\Pi_{N}$-formula, and 
$RM_{N}(\mathcal{X})$ the $\Pi_{N}$-recursively Mahlo operation for classes of transitive sets $\mathcal{X}$:
\begin{eqnarray*}
P\in RM_{N}(\mathcal{X}) & :\Leftrightarrow & \forall b\in P[P\models\Pi_{N}(b) \to \exists Q\in \mathcal{X}\cap P(Q\models\Pi_{N}(b))]
\\
&&\mbox{(read:} P \mbox{ is } \Pi_{N}\mbox{-reflecting on } \mathcal{X}\mbox{.)}
\end{eqnarray*}

The iteration of $RM_{N}$ along a definable relation $\prec$ is defined as follows.
\[
P\in RM_{N}(a;\prec) :\Leftrightarrow a\in P\in\bigcap\{ RM_{N}(RM_{N}(b;\prec)): b\in P\models b\prec a\}.
\]
Let $Ord\subset V$ denote the class of ordinals, $Ord^{\varepsilon}\subset V$ and $<^{\varepsilon}$ be $\Delta$-predicates such that
for any transitive and wellfounded model $V$ of $\mbox{{\sf KP}}\omega$,
$<^{\varepsilon}$ is a well ordering of type $\varepsilon_{I+1}$ on $Ord^{\varepsilon}$
for the order type $I$ of the class $Ord$ in $V$.
Specifically let us encode `ordinals' $\alpha<\varepsilon_{I+1}$ by codes $\lceil\alpha\rceil\in Ord^{\varepsilon}$ as follows.
$\lceil\alpha\rceil=\langle 0,\alpha\rangle$ for $\alpha\in Ord$, 
$\lceil I\rceil=\langle 1,0\rangle$,
$\lceil\omega^{\alpha}\rceil=\langle 2,\lceil\alpha\rceil\rangle$ for $\alpha>I$,
and
$\lceil \alpha\rceil=\langle 3,\lceil\alpha_{1}\rceil,\ldots,\lceil\alpha_{n}\rceil\rangle$
if $\alpha=\alpha_{1}+\cdots+\alpha_{n}>I$ with $\alpha_{1}\geq\cdots\geq\alpha_{n}$, $n>1$ and 
$\exists\beta_{i}(\alpha_{i}=\omega^{\beta_{i}})$ for each $\alpha_{i}$.
Then
$\lceil\omega_{n}(I+1)\rceil\in Ord^{\varepsilon}$ denotes the code of the `ordinal' 
$\omega_{n}(I+1)$.

$<^{\varepsilon}$ is assumed to be a canonical ordering such that 
$\mbox{{\sf KP}}\omega$ proves the fact that $<^{\varepsilon}$ is a linear ordering, and for any formula $\varphi$
and each $n<\omega$,
\begin{equation}\label{eq:trindveps}
\mbox{{\sf KP}}\omega\vdash\forall x(\forall y<^{\varepsilon}x\,\varphi(y)\to\varphi(x)) \to \forall x<^{\varepsilon}\lceil\omega_{n}(I+1)\rceil\varphi(x)
\end{equation}
For a definition of $\Delta$-predicates $Ord^{\varepsilon}$ and $<^{\varepsilon}$, 
and a proof of (\ref{eq:trindveps}), cf. \cite{liftupZF}.

\begin{theorem}\label{th:resolutionpi3}

For each $N\geq 2$,
{\rm {\sf KP}}$\Pi_{N+1}$ is $\Pi_{N+1}$-conservative over
the theory 
\[
\mbox{{\rm {\sf KP}}}\omega+
\{
V\in RM_{N}(\lceil\omega_{n}(I+1)\rceil;<^{\varepsilon}): 
n\in\omega\}
.\]
\end{theorem}

On the other side, we\cite{liftupZF} have lifted up the ordinal analysis of recursively inaccessible ordinals in \cite{Buchholz} 
to one of weakly inaccessible cardinals.
This paper aims to lift up \cite{Rathjen94} and \cite{LMPS} to the weak compactness.
\\

Let us mention the contents of this paper.
In the next section \ref{sect:ordinalnotation} 
iterated Skolem hulls $\mathcal{H}_{\alpha,n}(X)$ of sets $X$ of ordinals, ordinals $\Psi_{\kappa,n}\gamma$ for regular ordinals 
$\kappa\,(\mathcal{K}<\kappa\leq I)$,
and classes $Mh_{n}^{\alpha}[\Theta]$ are defined for finite sets $\Theta$ of ordinals.
It is shown that for each $n,m<\omega$, $(\mathcal{K} \mbox{ {\rm is a weakly compact cardinal}}) \to \mathcal{K}\in Mh_{n}^{\omega_{m}(I+1)}$ in
${\sf ZF}+(V=L)$.
In the third section \ref{sect:Ztheory} we introduce a theory for weakly compact cardinals,
which are equivalent to ${\sf ZF}+(V=L)+(\mathcal{K} \mbox{ {\rm is a weakly compact cardinal}}) $.

In the section \ref{sect:controlledOme}
cut inferences are eliminated from operator controlled derivations of $\Sigma^{1}_{2}$-sentences $\varphi^{V_{\mathcal{K}}}$ over $\mathcal{K}$, and $\varphi^{V_{\mathcal{K}}}$ is shown to be true.
Everything up to this is seen to be formalizable in 
${\sf ZF}+(V=L)+\{\theta_{n}(\mathcal{K}) : n\in\omega \}$.
Hence the Theorem \ref{th:mainthK} follows in the final section \ref{sect:prfmainthm}.

\section{Ordinals for weakly compact cardinals}\label{sect:ordinalnotation}
In this section
iterated Skolem hulls $\mathcal{H}_{\alpha,n}(X)$ of sets $X$ of ordinals, ordinals $\Psi_{\kappa,n}\gamma$ for regular ordinals 
$\kappa\,(\mathcal{K}<\kappa\leq I)$,
and classes $Mh_{n}^{\alpha}[\Theta]$ are defined for finite sets $\Theta$ of ordinals.
It is shown that for each $n,m<\omega$, $\mathcal{K}\in Mh_{n}^{\omega_{m}(I+1)}$ in
${\sf ZF}+(V=L)$ assuming $\mathcal{K}$ is a weakly compact cardinal.

Let $Ord^{\varepsilon}$ and $<^{\varepsilon}$ are $\Delta$-predicates 
as described before Theorem \ref{th:resolutionpi3}.
In the definition of $Ord^{\varepsilon}$ and $<^{\varepsilon}$,
$I$ with its code $\lceil I\rceil=\langle 1,0\rangle$ 
is {\it intended\/} to denote  the least weakly inaccessible cardinal
above the least weakly compact cardinal $\mathcal{K}$,
though we {\it do not assume\/} the existence of weakly inaccessible cardinals above $\mathcal{K}$
anywhere in this paper.
We are working in ${\sf ZF}+(V=L)$ assuming $\mathcal{K}$ is a weakly compact cardinal.

$Reg$ denotes the set of uncountable regular ordinals 
above $\mathcal{K}$, while
$R:=Reg\cap\{\rho:\mathcal{K}<\rho<I\}$ and
$R^{+}:=R\cup\{I\}$.
$\kappa,\lambda,\rho,\pi$ denote elements of $R$.
$\kappa^{+}$ denotes the least regular ordinal above $\kappa$.
$\Theta$ denotes finite sets of ordinals$\leq\mathcal{K}$.
$\Theta\subset_{fin} X$ iff $\Theta$ is a finite subset of $X$.
$Ord$ denotes the class of ordinals less than $I$, 
while $Ord^{\varepsilon}$ the class of codes of ordinals less than the next epsilon number $\varepsilon_{I+1}$ to $I$.

For admissible ordinals $\sigma$ and $X\subset L_{\sigma}$, $\mbox{{\rm Hull}}_{\Sigma_{n}}^{\sigma}(X)$ denotes  the $\Sigma_{n}$-Skolem hull of $X$ over $L_{\sigma}$, cf. \cite{liftupZF}.
$F(y)=F^{\Sigma_{n}}(y;\sigma,X)$ denotes the Mostowski collapsing
$F: \mbox{{\rm Hull}}_{\Sigma_{n}}^{\sigma}(X)\leftrightarrow L_{\gamma}$
of $\mbox{{\rm Hull}}_{\Sigma_{n}}^{\sigma}(X)$ for a $\gamma$.
Let $F^{\Sigma_{n}}(\sigma;\sigma,X):=\gamma$.
When $\sigma=I$, we write $F^{\Sigma_{n}}_{X}(y)$ for $F^{\Sigma_{n}}(y;I,X)$.

{\it In what follows\/} $n\geq 1$ {\it denotes a fixed positive integer\/}.
\\

$Code^{\varepsilon}$ denotes the union of codes $Ord^{\varepsilon}$ of ordinals$<\varepsilon_{I+1}$,
and codes $L_{I}:=\{\langle 0,x\rangle: x\in L\}$ of sets $x$ in the universe $L$.

For $\alpha,\beta\in Ord^{\varepsilon}$,
$\alpha\oplus\beta, \tilde{\omega}^{\alpha}\in Ord^{\varepsilon}$ 
denotes the codes of the sum and exponentiation, resp.

Let
\[
I:=\langle 1,0\rangle, \: \omega_{n}(I+1):=\tilde{\omega}_{n}( \langle 3,\langle 1,0\rangle, \langle 0,1\rangle\rangle), 
\mbox{ and } L_{I}:=\{\langle 0,x\rangle: x\in L\}
\]
and for codes $X,Y\in Code^{\varepsilon}$ 
\[
X\subset^{\varepsilon}Y:\Leftrightarrow \forall x\in^{\varepsilon} X(x\in^{\varepsilon}Y)
.\]

For simplicity let us identify the code $x\in Code^{\varepsilon}$ with
the `set' coded by $x$,
and $\in^{\varepsilon}$ [$<^{\varepsilon}$] is denoted by $\in$ [$<$], resp. when no confusion likely occurs.
For example, the code $\langle  0,x\rangle$ is identified with the set $\{\langle  0,y\rangle: y\in x\}$ of codes.

Define simultaneously 
 the classes $\mathcal{H}_{\alpha,n}(X)\subset L_{I}\cup\{x\in Ord^{\varepsilon}: x<^{\varepsilon}\omega_{n+1}(I+1)\}$, 
and the ordinals $\Psi_{\kappa,n} \alpha\,(\kappa\in R^{+})$
for $\alpha<^{\varepsilon}\omega_{n+1}(I+1)$ and sets $X\subset L_{I}$ as follows.
We see that $\mathcal{H}_{\alpha,n}(X)$ and $\Psi_{\kappa,n} \alpha$ are (first-order) definable as a fixed point in ${\sf ZF}+(V=L)$ cf. Proposition \ref{prp:definability}.

$\mathcal{H}_{\alpha,n}$ is an operator in the sense defined below.

\begin{definition}\label{df:operator}
{\rm By an} {\it operator\/} {\rm we mean a map} $\mathcal{H}$, 
$\mathcal{H}:\mathcal{P}(L_{I})\to\mathcal{P}(L_{I}\cup\{x\in Ord^{\varepsilon}: x<^{\varepsilon}\omega_{n+1}(I+1)\})${\rm , such that}
\benu
\item
$\forall X\subset L_{I}[X\subset\mathcal{H}(X)]$.

\item
$\forall X,Y\subset L_{I}[Y\subset\mathcal{H}(X) \Rightarrow \mathcal{H}(Y)\subset\mathcal{H}(X)]$.
\eenu

{\rm For an operator} $\mathcal{H}$ {\rm and} $\Theta,\Lambda\subset L_{I}$,
$\mathcal{H}[\Theta](X):=\mathcal{H}(X\cup\Theta)${\rm , and}
$\mathcal{H}[\Theta][\Lambda]:=(\mathcal{H}[\Theta])[\Lambda]${\rm , i.e.,}
$\mathcal{H}[\Theta][\Lambda](X)=\mathcal{H}(X\cup\Theta\cup\Lambda)$.
\end{definition}
Obviously $\mathcal{H}[\Theta]$ is an operator.

\begin{definition}\label{df:Cpsiregularsm}

$\mathcal{H}_{\alpha,n}(X)$ {\rm is a} 
{\rm Skolem hull of} $\{\langle 0,0\rangle,\mathcal{K},I\}\cup X$ {\rm under the functions} 
$\oplus,
 \alpha\mapsto\tilde{\omega}^{\alpha}, 
 \kappa\mapsto\kappa^{+}\, (\kappa\in R),
 \Psi_{\kappa,n}\!\upharpoonright\! \alpha\,(\kappa\in R^{+})${\rm , the Skolem hullings:}
\[
X \mapsto  \mbox{{\rm Hull}}^{I}_{\Sigma_{n}}(X\cap I)
\]
{\rm and the Mostowski collapsing functions}
\[
x=\Psi_{\kappa,n}\gamma\mapsto F^{\Sigma_{1}}_{x\cup\{\kappa\}}\,(\kappa\in R)
\]
{\rm and}
\[
x=\Psi_{I,n}\gamma\mapsto F^{\Sigma_{n}}_{x}
\]

\benu

\item
{\rm (Inductive definition of} $\mathcal{H}_{\alpha,n}(X)${\rm ).}
 \benu
\item
$\{\langle 0,0\rangle,\mathcal{K},I\}\cup X\subset\mathcal{H}_{\alpha,n}(X)$.

\item
 $x, y \in \mathcal{H}_{\alpha}(X) \Rightarrow x\oplus y,\tilde{\omega}^{x}\in \mathcal{H}_{\alpha,n}(X)$.
 
\item
$\kappa\in\mathcal{H}_{\alpha,n}(X)\cap (\{\mathcal{K}\}\cup R) \Rightarrow \kappa^{+}\in\mathcal{H}_{\alpha,n}(X)$.

\item
$
\gamma\in \mathcal{H}_{\alpha,n}(X)\cap\alpha
\Rightarrow 
\Psi_{I,n}\gamma\in\mathcal{H}_{\alpha,n}(X)
$.

\item

{\rm If}
$\kappa\in\mathcal{H}_{\alpha,n}(X)\cap R$,
$\gamma\in \mathcal{H}_{\alpha,n}(X)\cap\alpha$ 
{\rm and} $\kappa\in\mathcal{H}_{\gamma,n}(\kappa)$,
{\rm then}
$\Psi_{\kappa,n}\gamma\in\mathcal{H}_{\alpha,n}(X)$.

\item

\[
\mbox{{\rm Hull}}^{I}_{\Sigma_{n}}(\mathcal{H}_{\alpha,n}(X)\cap L_{I})\cap Code^{\varepsilon} 
\subset \mathcal{H}_{\alpha,n}(X)
.\]
{\rm Namely for any} $\Sigma_{n}${\rm -formula} $\varphi[x,\vec{y}\,]$ {\rm in the language} $\{\in\}$
{\rm and parameters} $\vec{a}\subset \mathcal{H}_{\alpha,n}(X)\cap L_{I}${\rm , if} $b\in L_{I}$,
$(L_{I},\in^{\varepsilon})\models\varphi[b,\vec{a}\,]$ {\rm and} 
$(L_{I},\in^{\varepsilon})\models\exists!x\,\varphi[x,\vec{a}\,]$,
{\rm then} $b\in\mathcal{H}_{\alpha,n}(X)$.

\item
{\rm If} $\kappa\in\mathcal{H}_{\alpha,n}(X)\cap R$, 
$\gamma\in \mathcal{H}_{\alpha,n}(X)\cap\alpha$, $x=\Psi_{\kappa,n}\gamma\in\mathcal{H}_{\alpha,n}(X)$, 
$\kappa\in\mathcal{H}_{\gamma,n}(\kappa)$
{\rm and}
$\delta\in (\mbox{{\rm Hull}}^{I}_{\Sigma_{1}}(x\cup\{\kappa\})\cup\{I\})\cap\mathcal{H}_{\alpha,n}(X)$,
{\rm then}
$F^{\Sigma_{1}}_{x\cup\{\kappa\}}(\delta)\in\mathcal{H}_{\alpha,n}(X)$.

\item
{\rm If} 
$\gamma\in \mathcal{H}_{\alpha,n}(X)\cap\alpha$, $x=\Psi_{I,n}\gamma\in\mathcal{H}_{\alpha,n}(X)$, 
{\rm and}
$\delta\in (\mbox{{\rm Hull}}^{I}_{\Sigma_{n}}(x)\cup\{I\})\cap\mathcal{H}_{\alpha,n}(X)$,
{\rm then}
$F^{\Sigma_{n}}_{x}(\delta)\in\mathcal{H}_{\alpha,n}(X)$.

 \eenu

\item
{\rm (Definition of} $\Psi_{\kappa,n}\alpha${\rm ).}

{\rm Assume} $\kappa\in R^{+}$ {\rm and} $\kappa\in\mathcal{H}_{\alpha,n}(\kappa)${\rm . Then}
\[
\Psi_{\kappa,n}\alpha:=
\min_{^{\varepsilon}}\{\beta<^{\varepsilon}\kappa : \kappa\in \mathcal{H}_{\alpha,n}(\beta),\, \mathcal{H}_{\alpha,n}(\beta)\cap \kappa \subset^{\varepsilon}\beta\} 
.\]

\eenu

\end{definition}

Definition \ref{df:Cpsiregularsm} is essentially the same as in \cite{liftupZF}.

The classes $Mh^{\alpha}_{n}[\Theta]$ are defined for $n<\omega$, 
$\alpha<\varepsilon_{I+1}$, and $\Theta\subset_{fin}(\mathcal{K}+1)$.

\begin{definition}\label{df:dfMh}
{\rm (}$Mh^{\alpha}_{n}[\Theta]${\rm )}

{\rm Let} $\Theta\subset_{fin}(\mathcal{K}+1)$ {\rm and} $\mathcal{K}\geq\pi\in Reg ${\rm . Then}
\begin{eqnarray}
\pi\in Mh^{\alpha}_{n}[\Theta] & :\Leftrightarrow &
\mathcal{H}_{\alpha,n}(\pi)\cap\mathcal{K}\subset^{\varepsilon}\pi
\,\&\,
 \alpha\in \mathcal{H}_{\alpha,n}[\Theta](\pi)
\nonumber
\\
& \,\&\,  &
\forall\xi \in \mathcal{H}_{\xi,n}[\Theta\cup\{\pi\}](\pi)\cap\alpha[
\pi\in M(Mh_{n}^{\xi}[\Theta\cup\{\pi\}])]
\label{eq:dfMh}
\end{eqnarray}
{\rm where}
$\forall\xi \in  \mathcal{H}_{\xi,n}[\Theta\cup\{\pi\}](\pi)\cap\alpha[\cdots]$ {\rm is a short hand for}
$\forall\xi <^{\varepsilon}\alpha[\xi\in  \mathcal{H}_{\xi,n}[\Theta\cup\{\pi\}](\pi)\cap\alpha\to \cdots]$.

\[
Mh_{n}^{\alpha}  :=  Mh_{n}^{\alpha}[\{\mathcal{K}\}]=Mh_{n}^{\alpha}[\emptyset]
.\]
\end{definition}

The following Propositions \ref{prp:definability} and \ref{prp:clshull.-1} are easy to see.

\begin{proposition}\label{prp:definability}
Each of 
$x=\mathcal{H}_{\alpha,n}(\beta)\, (\alpha\in Ord^{\varepsilon},\beta<^{\varepsilon}I)$,
$\beta=\Psi_{\kappa,n}\alpha\,(\kappa\in R^{+})$ and
$x=Mh_{n}^{\alpha}[\Theta]$
is a $\Sigma_{n+1}$-predicate as fixed points in ${\sf ZF}+(V=L)$. 
\end{proposition}

\begin{proposition}\label{prp:clshull.-1}
$(\alpha,y)\mapsto \mathcal{H}_{\alpha,n}[\Theta](y)$ is weakly monotonic in the sense that
\[
\alpha\leq^{\varepsilon}\alpha^{\prime}\land y\subset y^{\prime} \land
x=\mathcal{H}_{\alpha,n}[\Theta](y)\land x^{\prime}=\mathcal{H}_{\alpha^{\prime},n}[\Theta](y^{\prime})\to x\subset x^{\prime}
.\]

Also $(\alpha,y)\mapsto \mathcal{H}_{\alpha,n}[\Theta](y)$ is continuous in the sense that
 if $\alpha=\sup_{i\in I}\alpha_{i}$ is a limit ordinal with an increasing sequence $\{\alpha_{i}\}_{i\in I}$ and
 $y=\bigcup_{j\in J}y_{j}$ with a directed system $\{y_{j}\}_{j\in J}$, then
\[
x=\mathcal{H}_{\alpha,n}[\Theta](\beta)\land \forall i\in I\forall j\in J(x_{i,j}=\mathcal{H}_{\alpha_{i},n}[\Theta](y_{j}))
\to  x=\bigcup_{i\in I,j\in J}x_{i,j}
.\]
\end{proposition}

Let $A_{n}(\alpha)$ denote the conjunction of 
$\forall\beta<^{\varepsilon}I\exists ! x[x=\mathcal{H}_{\alpha,n}(\beta)]$,
$\forall\kappa\in R^{+}\forall x[\kappa\in x=\mathcal{H}_{\alpha,n}(\kappa) \to \exists!\beta(\beta=\Psi_{\kappa,n}\alpha)]$ and
$\forall\Theta\subset_{fin}(\mathcal{K}+1)\exists ! x[x=Mh_{n}^{\alpha}[\Theta]]$.
\\

The $\Sigma_{n+1}$-formula $\theta_{n}(x)$ in Theorem \ref{th:mainthK}
is defined to be
\[
\theta_{n}(x) :\equiv \exists y[y=Mh_{n}^{\omega_{n}(I+1)} \land x\in y] 
.\]
The following Lemma \ref{lem:welldefinedness}.\ref{lem:welldefinedness.1} shows Theorem \ref{th:mainthK}.\ref{th:mainthK1}.

$card(x)$ denotes the cardinality of sets $x$.

\begin{lemma}\label{lem:welldefinedness}
For each $n,m<\omega$, ${\sf ZF}+(V=L)$ proves the followings.
\benu

\item\label{lem:welldefinedness.0}
$y=\mathcal{H}_{\alpha,n}(x) \to card(y)\leq\max\{card(x),\aleph_{0}\}$.

\item\label{lem:welldefinedness.2}
$\forall\alpha<^{\varepsilon}\omega_{m}(I+1)\, A_{n}(\alpha)
$.

\item\label{lem:welldefinedness.1}
If $\mathcal{K}$ is weakly compact and $\Theta\subset_{fin}(\mathcal{K}+1)$, then
$\mathcal{K}\in Mh_{n}^{\omega_{m}(I+1)}[\Theta]\cap M(Mh_{n}^{\omega_{m}(I+1)}[\Theta])
$.
\eenu

\end{lemma}
\bprf

\ref{lem:welldefinedness}.\ref{lem:welldefinedness.2}.
We show that $A_{n}(\alpha)$ is progressive, i.e.,
$\forall\alpha<^{\varepsilon}\omega_{m}(I+1)[\forall\gamma<^{\varepsilon}\alpha\, A_{n}(\gamma) \to A_{n}(\alpha)]$.

Assume $\forall\gamma<^{\varepsilon}\alpha\, A_{n}(\gamma)$ and $\alpha<^{\varepsilon}\omega_{m}(I+1)$.
$\forall\beta<^{\varepsilon}I\exists ! x[x=\mathcal{H}_{\alpha,n}(\beta)]$
 follows from IH and the Replacement.

Next assume $\kappa\in R^{+}$ and $\kappa\in\mathcal{H}_{\alpha,n}(\kappa)$.
Then $\exists!\beta(\beta=\Psi_{\kappa,n}\alpha)$ follows from the regularity of $\kappa$ and Proposition \ref{prp:clshull.-1}.

$\exists ! x[x=Mh_{n}^{\alpha}[\Theta]]$ is easily seen from IH.
\\

\noindent
\ref{lem:welldefinedness}.\ref{lem:welldefinedness.1}.
Suppose $\mathcal{K}$ is $\Pi^{1}_{1}$-indescribable.
We show 
\[
B_{n}(\alpha):\Leftrightarrow \forall\Theta\subset_{fin}(\mathcal{K}+1)[\alpha\in \mathcal{H}_{\alpha,n}[\Theta](\mathcal{K})
\to \mathcal{K}\in Mh_{n}^{\alpha}[\Theta] \cap M(Mh_{n}^{\alpha}[\Theta] )]
\]
is progressive in $\alpha$.
 
Suppose $\forall\xi<^{\varepsilon}\alpha\, B_{n}(\xi)$, $\Theta\subset_{fin}(\mathcal{K}+1)$ and $\alpha\in \mathcal{H}_{\alpha,n}[\Theta](\mathcal{K})$.
We have to show that $Mh_{n}^{\alpha}[\Theta]$ meets every club subset $C_{0}$ of $\mathcal{K}$.
$\mathcal{K}\in Mh_{n}^{\alpha}[\Theta]$ follows from $\mathcal{K}\in M(Mh_{n}^{\alpha}[\Theta])$, cf. Proposition \ref{prp:clshull}.\ref{prp:Mh2}.
We can assume that 
$\forall\pi\in C_{0}[ (\mathcal{H}_{\alpha,n}(\pi)\cap\mathcal{K}\subset\pi) \land(\alpha\in\mathcal{H}_{\alpha,n}[\Theta](\pi))]$
since both of $\{\pi<\mathcal{K}:\mathcal{H}_{\alpha,n}(\pi)\cap\mathcal{K}\subset\pi\}$ and $\{\pi<\mathcal{K}: \alpha\in\mathcal{H}_{\alpha,n}[\Theta](\pi)\}$
are club in $\mathcal{K}$.

Since $\forall\pi\leq\mathcal{K}[card(\mathcal{H}_{\alpha,n}[\Theta\cup\{\pi\}](\pi))\leq\pi]$, 
pick an injection $f:\mathcal{H}_{\alpha,n}[\Theta\cup\{\mathcal{K}\}](\mathcal{K})\to \mathcal{K}$ so that
$f"\mathcal{H}_{\alpha,n}[\Theta\cup\{\pi\}](\pi)\subset\pi$ for any weakly inaccessibles $\pi\leq\mathcal{K}$.

Let $R_{0}=\{f(\alpha)\}$, $R_{1}=C_{0}$, 
$R_{2}=\{f(\xi) :\xi\in \mathcal{H}_{\xi,n}[\Theta](\mathcal{K})\cap\alpha\}$,
$R_{3}=\bigcup\{(Mh_{n}^{\xi}[\Theta\cup\{\pi\}]\cap\mathcal{K})\times\{f(\pi)\}\times\{f(\xi)\} :\xi\in \mathcal{H}_{\xi,n}[\Theta](\mathcal{K})\cap\alpha, \pi\leq\mathcal{K}\}$,
and $R_{4}=\{(f(\beta),f(\gamma)): \{\beta,\gamma\}\subset\mathcal{H}_{\alpha,n}[\Theta\cup\{\mathcal{K}\}](\mathcal{K}), \beta<\gamma\}$.

By IH we have $\forall\xi\in \mathcal{H}_{\xi,n}[\Theta](\mathcal{K})\cap\alpha[\mathcal{K}\in M(Mh_{n}^{\xi}[\Theta])]$.
Hence $\langle V_{\mathcal{K}},\in, R_{i}\rangle_{i\leq 4}$
enjoys a $\Pi^{1}_{1}$-sentence saying that
$\mathcal{K}$ is weakly inaccessible, $R_{0}\neq\emptyset$, $R_{1}$ is a club subset of $\mathcal{K}$ and
\[
\varphi:\Leftrightarrow \forall C\mbox{:club } \forall x, y [
R_{2}(x) \land \theta(R_{4},y)
\to C\cap \{a: R_{3}(a,y,x)\}\neq\emptyset]
\]
where $\theta(R_{4},y)$ is a $\Sigma^{1}_{1}$-formula such that for any $\pi\leq\mathcal{K}$
\[
V_{\pi}\models \theta(R_{4},y) \Leftrightarrow y=f(\pi)
\]
Namely $\theta(R_{4},y)$ says that there exists a function $G$ on the class $Ord$ of ordinals such that 
$
\forall \beta,\gamma\in Ord[(\beta<\gamma \leftrightarrow R_{4}(G(\beta),G(\gamma))\land(G(\beta)<y)]
$ and
$\forall z(R_{4}(z,y) \to \exists\beta\in Ord(G(\beta)=z))
$.

By the $\Pi^{1}_{1}$-indescribability of $\mathcal{K}$, pick a $\pi<\mathcal{K}$ such that
$\langle V_{\pi},\in, R_{i}\cap V_{\pi}\rangle_{i\leq 4}$ enjoys the $\Pi^{1}_{1}$-sentence.

We claim $\pi\in C_{0}\cap Mh_{n}^{\alpha}[\Theta]$.
$\pi$ is weakly inaccessible, $f(\alpha)\in V_{\pi}$ and $C_{0}$ is club in $\pi$, and hence $\pi\in C_{0}$.
It remains to see $\forall\xi\in \mathcal{H}_{\xi,n}[\Theta\cup\{\pi\}](\pi)\cap\alpha[\pi\in M(Mh_{n}^{\xi}[\Theta\cup\{\pi\}])]$.
This follows from the fact that
$\varphi$ holds in $\langle V_{\pi},\in, R_{i}\cap V_{\pi}\rangle_{i\leq 4}$, and
$\forall\xi\in \mathcal{H}_{\xi,n}[\Theta\cup\{\pi\}](\pi)\cap\alpha(f(\xi)\in V_{\pi})$
by $f" \mathcal{H}_{\alpha,n}[\Theta\cup\{\pi\}](\pi)\subset\pi$ and $\mathcal{H}_{\xi,n}[\Theta\cup\{\pi\}](\pi)\subset\mathcal{H}_{\xi,n}[\Theta](\mathcal{K})$.

Thus $\mathcal{K}\in M(Mh_{n}^{\alpha}[\Theta])$.
\eprf

\begin{definition}\label{df:Hn}
 $\mathcal{H}(n)$ {\rm denotes a subset of} $\mathcal{H}_{\omega_{n}(I+1),n}(\emptyset)$ {\rm such that
every ordinal is hereditarily less than} $\omega_{n}(I+1)$.

{\rm This means} $\alpha\in\mathcal{H}(n)\Rightarrow \alpha<\omega_{n}(I+1)${\rm , etc.}
\end{definition}

\bcor\label{cor:existenceord}
For each $n<\omega$, $\mathcal{H}(n)$ is well-defined in ${\sf ZF}+(V=L)$.
\ecor

Let us see some elementary facts.

\begin{proposition}\label{prp:clshull}

\benu

\item\label{prp:Mh1}
$\alpha\in\mathcal{H}_{\alpha,n}[\Theta](\pi) \,\&\, \pi\in Mh_{n}^{\alpha}[\Theta\cup\{\rho\}]
\Rightarrow \pi\in Mh_{n}^{\alpha}[\Theta]$.

\item\label{prp:Mh2}
$\pi\in M(Mh_{n}^{\alpha}[\Theta\cup\{\pi\}])\Rightarrow \pi\in Mh_{n}^{\alpha}[\Theta\cup\{\pi\}]$.

\item\label{prp:Mh3}
$\pi\in Mh_{n}^{\alpha}[\Theta] \,\&\, \xi\in \mathcal{H}_{\xi,n}[\Theta\cup\{\pi\}](\pi)\cap\alpha \Rightarrow \pi\in Mh_{n}^{\xi}[\Theta\cup\{\pi\}]$, and
$\pi\in Mh_{n}^{\alpha}[\Theta] \,\&\, \xi\in \mathcal{H}_{\xi,n}[\Theta](\pi)\cap\alpha \Rightarrow \pi\in Mh_{n}^{\xi}[\Theta]$.

\eenu

\end{proposition}
\bprf

\ref{prp:clshull}.\ref{prp:Mh2}. This is seen from Proposition \ref{prp:clshull}.\ref{prp:Mh1}.
\\

\noindent
\ref{prp:clshull}.\ref{prp:Mh3}. This is seen from Proposition \ref{prp:clshull}.\ref{prp:Mh2}.

\eprf

\subsection{Greatly Mahlo cardinals}

Let us compare the class $Mh^{\alpha}_{n}[\Theta]$ with Rathjen's class $M^{\alpha}$ in \cite{Rathjen94}.
The difference lies in augmenting finite sets $\Theta$ of ordinals, which are given in advance.
Moreover the finite set grows when we step down to previously defined classes, cf. (\ref{eq:dfMh}).
For example if 
an ordinal $\xi<\alpha$ is $\Sigma_{1}$-definable from $\{\pi,\pi^{+}\}$,
then $\xi \in \mathcal{H}_{\xi,n}[\Theta\cup\{\pi\}](\pi)$ for $n\geq 1$.
Hence $Mh_{n}^{\xi}[\Theta\cup\{\pi\}]$ is stationary in $\pi$ for such an ordinal $\xi<\alpha$
if $\pi\in Mh^{\alpha}_{n}[\Theta]$.
Cf. {\bf Case 2} in the proof of Lemma \ref{lem:CollapsingthmKR} below.
 
This yields that 
any $\sigma$ with $\sigma\in   Mh_{n}^{\sigma^{+}}$ is a greatly Mahlo cardinal in the sense of 
Baumgartner-Taylor-Wagon\cite{Baumgartner}.
Moreover if $\mathcal{K}\in Mh_{n}^{\mathcal{K}+1}$, then 
the class of the greatly Mahlo cardinals below $\mathcal{K}$ is stationary in $\mathcal{K}$
as seen in Proposition \ref{prp:Baumgarter}.

$M^{\alpha}\, (\alpha<\mathcal{K}^{+})$ denotes the set of $\alpha$-weakly Mahlo cardinals defined as follows.
$M^{0}:=Reg \cap\mathcal{K}$,
$M^{\alpha+1}=M(M^{\alpha})$,
$M^{\lambda}=\bigcap\{M(M^{\alpha}): \alpha<\lambda\}$ for limit ordinals $\lambda$ with $cf(\lambda)<\mathcal{K}$,
and
$M^{\lambda}:=\triangle\{M(M^{\lambda_{i}}):i<\mathcal{K}\}$ for limit ordinals $\lambda$ with $cf(\lambda)=\mathcal{K}$,
where $\sup_{i<\mathcal{K}}\lambda_{i}=\lambda$ and the sequence $\{\lambda_{i}\}_{i<\mathcal{K}}$ is chosen so that
it is the $<_{L}$-minimal such sequence.

In the last case for $\pi<\mathcal{K}$,
$\pi\in M^{\lambda}\Leftrightarrow \forall i<\pi(\pi\in M(M^{\lambda_{i}}))$.

\begin{proposition}\label{prp:Baumgarter}
For $n\geq 1$ and $\sigma\leq\mathcal{K}$, the followings are provable in ${\sf ZF}+(V=L)$.
\benu

\item\label{prp:Baumgarter.1}
If $\sigma\in\Theta$,
$\pi\in Mh_{n}^{\alpha}[\Theta]\cap\sigma$, and $\alpha\in\mbox{{\rm Hull}}^{I}_{\Sigma_{1}}(\{\sigma,\sigma^{+}\}\cup\pi)\cap\sigma^{+}$, then
$\pi\in M^{\alpha}$.

\item\label{prp:Baumgarter.2}
$
\sigma\in Mh_{n}^{\sigma^{+}}[\Theta] \to \forall\alpha<\sigma^{+}(\sigma\in M(M^{\alpha}))
$.

\item\label{prp:Baumgarter.3}
The class of the greatly Mahlo cardinals below $\mathcal{K}$ is stationary in $\mathcal{K}$
if $\mathcal{K}\in Mh_{n}^{\mathcal{K}+1}$.
\eenu

\end{proposition}
\bprf

\ref{prp:Baumgarter}.\ref{prp:Baumgarter.1} by induction on $\alpha<\sigma^{+}$.
Suppose $\sigma\in\Theta$, 
$\pi\in Mh_{n}^{\alpha}[\Theta]\cap\sigma$ and $\alpha\in\mbox{Hull}^{I}_{\Sigma_{1}}(\{\sigma,\sigma^{+}\}\cup\pi)\cap\sigma^{+}$.

First consider the case when $cf(\alpha)=\sigma$, and let $\{\alpha_{i}\}_{i<\sigma}$ be the $<_{L}$-minimal
sequence such that $\sup_{i<\sigma}\alpha_{i}=\alpha$.
Then $\{\alpha_{i}\}_{i<\sigma}\in\mbox{Hull}^{I}_{\Sigma_{1}}(\{\alpha,\sigma\})\subset\mbox{Hull}^{I}_{\Sigma_{1}}(\{\sigma,\sigma^{+}\}\cup\pi)$.
For $i<\pi$,
$\alpha_{i}\in \mbox{Hull}^{I}_{\Sigma_{1}}(\{\sigma,\sigma^{+}\}\cup\pi)\cap\alpha\subset\mathcal{H}_{\alpha_{i},n}[\Theta\cup\{\pi\}](\pi)\cap\alpha$
by $\sigma\in\Theta$.
$\pi\in Mh_{n}^{\alpha}[\Theta]$ yields $\pi\in M(Mh_{n}^{\alpha_{i}}[\Theta\cup\{\pi\}])$.
Now for a club subset $C$ in $\pi$, pick a $\rho<\pi$ such that $\rho\in C\cap Mh_{n}^{\alpha_{i}}[\Theta\cup\{\pi\}]$.
We can assume that $\alpha_{i}\in\mbox{Hull}^{I}_{\Sigma_{1}}(\{\sigma,\sigma^{+}\}\cup\rho)$ 
by $\alpha_{i}\in \mbox{Hull}^{I}_{\Sigma_{1}}(\{\sigma,\sigma^{+}\}\cup\pi)$.
Thus IH yields $\rho\in M^{\alpha_{i}}$, and hence $\pi\in M(M^{\alpha_{i}})$ for any $i<\pi$.

Second consider the case when $cf(\alpha)<\sigma$.
Then $cf(\alpha)\in \mbox{Hull}^{I}_{\Sigma_{1}}(\{\alpha\})\cap\sigma\subset\mbox{Hull}^{I}_{\Sigma_{1}}(\{\sigma,\sigma^{+}\}\cup\pi)\cap\sigma\subset\mathcal{H}_{\alpha,n}[\{\sigma\}](\pi)\cap\sigma\subset\pi$ by $\pi\in Mh_{n}^{\alpha}[\Theta]$ and $\sigma\in\Theta$.
Thus $cf(\alpha)<\pi$.
Pick a cofinal sequence $\{\alpha_{i}\}_{i<cf(\alpha)}\in\mbox{Hull}^{I}_{\Sigma_{1}}(\{\sigma,\sigma^{+}\}\cup\pi)$.
Then for any $i<cf(\alpha)<\pi$
we have $\alpha_{i}\in \mbox{Hull}^{I}_{\Sigma_{1}}(\{\sigma,\sigma^{+}\}\cup\pi)\cap\alpha$,
and hence $\pi\in M(Mh_{n}^{\alpha_{i}}[\Theta\cup\{\pi\}])$.
As in the first case we see that $\pi\in M(M^{\alpha_{i}})$ for any $i<cf(\alpha)$.

Finally let $\alpha=\beta+1$.
Then $\beta\in\mbox{Hull}^{I}_{\Sigma_{1}}(\{\sigma,\sigma^{+}\}\cup\pi)$ together with IH yields $\pi\in M(M^{\beta})$.
\\

\noindent
\ref{prp:Baumgarter}.\ref{prp:Baumgarter.2}.
Suppose $\sigma\in Mh_{n}^{\sigma^{+}}[\Theta]$ and $\exists\alpha<\sigma^{+}(\sigma\not\in M(M^{\alpha}))$.
Let $\alpha<\sigma^{+}$ be the minimal ordinal such that $\sigma\not\in M(M^{\alpha})$, and $C$ be a club subset of $\sigma$
such that $C\cap M^{\alpha}=\emptyset$.
Then $\alpha\in\mbox{Hull}^{I}_{\Sigma_{1}}(\{\sigma,\sigma^{+}\})\cap\sigma^{+}\subset\mathcal{H}_{\alpha,n}[\Theta\cup\{\sigma\}](\sigma)\cap\sigma^{+}$.
By $\sigma\in Mh_{n}^{\sigma^{+}}[\Theta]$ we have $\sigma\in M(Mh_{n}^{\alpha}[\Theta\cup\{\sigma\}])$.
Pick a $\pi\in C\cap Mh_{n}^{\alpha}[\Theta\cup\{\sigma\}]\cap\sigma$. 
Proposition \ref{prp:Baumgarter}.\ref{prp:Baumgarter.1}
yields $\pi\in M^{\alpha}$. A contradiction.
\\

\noindent
\ref{prp:Baumgarter}.\ref{prp:Baumgarter.3}.
If $\mathcal{K}\in Mh_{n}^{\mathcal{K}+1}$, then $\mathcal{K}\in M(Mh_{n}^{\mathcal{K}})$.
Let $\sigma\in Mh_{n}^{\mathcal{K}}\cap\mathcal{K}$. 
Then $\sigma^{+}\in\mathcal{H}_{\sigma^{+},n}[\{\sigma\}](\sigma)\cap\mathcal{K}$, and hence
$\sigma\in M(Mh_{n}^{\sigma^{+}}[\{\sigma\}])$.
Proposition \ref{prp:clshull}.\ref{prp:Mh2} yields $\sigma\in Mh_{n}^{\sigma^{+}}[\{\sigma\}]$.
From Proposition \ref{prp:Baumgarter}.\ref{prp:Baumgarter.2}
we see that $\sigma$ is greatly Mahlo.
\eprf

\section{A theory for weakly compact cardinals}\label{sect:Ztheory}

In this section the set theory ${\sf ZF}+(V=L)+(\mathcal{K}\mbox{ is weakly compact)}$
is paraphrased to another set theory $\mbox{{\rm T}}(\mathcal{K},I)$ as in \cite{liftupZF}.

Let $\mathcal{K}$ be the least weakly compact cardinal, and $I>\mathcal{K}$ the least weakly inaccessible cardinal above $\mathcal{K}$.
$\kappa,\lambda,\rho$ ranges over uncountable regular ordinals such that $\mathcal{K}<\kappa,\lambda,\rho< I$.

In the following Definition \ref{df:regext}, 
the predicate $P$ is intended to denote the relation
\[
P(\lambda,x,y)  \Leftrightarrow  x=F^{\Sigma_{1}}_{x\cup\{\lambda\}}(\lambda)
\,\&\,
y=F^{\Sigma_{1}}_{x\cup\{\lambda\}}(I):=rng(F^{\Sigma_{1}}_{x\cup\{\lambda\}})\cap Ord
\]
and the predicate $P_{I,n}(x)$ is intended to denote the relation
\[
P_{I,n}(x)  \Leftrightarrow  x=F^{\Sigma_{n}}_{x}(I)
.\]

\begin{definition}\label{df:stratified}
\benu

\item\label{df:stratified1}
{\rm Let} $\vec{X}=X_{0},\ldots,X_{n-1}$ {\rm be a list of unary predicates.}
{\rm A} {\it stratified formula with respect to the variables\/} $\vec{x}=x_{0},\ldots,x_{n-1}$
 {\rm is a formula} $\varphi[\vec{x}\,]$ {\rm in the language} $\{\in\}$ 
{\rm obtained from a (first-order) formula} $\varphi[\vec{X}\,]$ {\rm  in the language} $\{\in\}\cup\vec{X}$ 
{\rm by replacing
any atomic formula} $X_{i}(z)$ {\rm by} $z\in x_{i}$ {\rm for} $i<n$.

\item\label{df:stratified2}
{\rm For a formula} $\varphi$ {\rm and a set} $x$,
$\varphi^{x}$ {\rm denotes the result of restricting every unbounded quantifier}
$\exists z,\forall z$ {\rm in} $\varphi$ {\rm to} $\exists z\in x, \forall z\in x$.

\item
$\alpha\in Ord:\Leftrightarrow \forall x\in a\forall y\in x(y\in a) \land
\forall x,y\in a(x\in y\lor x=y\lor y\in x)${\rm , and by}
$\alpha<\beta$ {\rm we tacitly assume that} $\alpha,\beta$ {\rm are ordinals, i.e.,}
$\alpha<\beta:\Leftrightarrow\{\alpha,\beta\}\subset Ord\land \alpha\in\beta$.
\eenu

\end{definition}

\begin{definition}\label{df:regext}
$\mbox{{\rm T}}(\mathcal{K},I,n)$ {\rm denotes the set theory defined as follows.}
\benu

\item
{\rm Its language is} $\{\in, P,P_{I,n},Reg,\mathcal{K}\}$ {\rm for a ternary predicate} $P${\rm  , unary predicates} $P_{I,n}$
{\rm and} $Reg${\rm , and an individual constant} $\mathcal{K}$.

\item
{\rm Its axioms are obtained from those of Kripke-Platek set theory with the axiom of infinity} $\mbox{{\rm KP}}\omega$
{\rm in the expanded language,
the axiom of constructibility,}
$V=L$
{\rm together with the axiom schemata saying that}
 \benu
 \item
 {\rm the ordinals} $\kappa$ 
{\rm with} $Reg(\kappa)$ {\rm is an uncountable regular ordinal}$>\mathcal{K}$
$(Reg(\kappa) \to \mathcal{K}<\kappa\in Ord)$ {\rm and } $(Reg(\kappa) \to a\in Ord\cap\kappa \to \exists x, y\in Ord\cap\kappa[a<x\land P(\kappa,x,y)])$,
{\rm and the ordinal}  $x$ {\rm with}
$P(\kappa,x,y)$ {\rm is a critical point of the} $\Sigma_{1}$ {\rm elementary embedding from an}
$L_{y}\cong \mbox{{\rm Hull}}^{I}_{\Sigma_{1}}(x\cup\{\kappa\})$ {\rm to the universe}
$L_{I}$ {\rm (}$P(\kappa,x,y) \to \{x,y\}\subset Ord \land  x<y<\kappa \land Reg(\kappa)$ {\rm and}  
$P(\kappa,x,y) \to a\in Ord\cap x \to \varphi[\kappa,a] \to \varphi^{y}[x,a]$
{\rm for any} $\Sigma_{1}${\rm -formula} $\varphi$ {\rm in the language} $\{\in\}${\rm ),}
\item
{\rm there are cofinally many regular ordinals (}$\forall x\in Ord \exists y[x\geq\mathcal{K} \to y>x\land Reg(y)]${\rm ),}
 \item
 {\rm the ordinal} $x$ {\rm with}
$P_{I,n}(x)$ {\rm is a critical point of the} $\Sigma_{n}$ {\rm elementary embedding from}
$L_{x}\cong \mbox{{\rm Hull}}^{I}_{\Sigma_{n}}(x)$ {\rm to the universe}
$L_{I}$ {\rm (}$P_{I,n}(x)\to x\in Ord$ {\rm and} $P_{I,n}(x) \to a\in Ord\cap x \to \varphi[a] \to \varphi^{x}[a]$
{\rm for any} $\Sigma_{n}${\rm -formula} $\varphi$ {\rm in the language} $\{\in\}${\rm ),}
{\rm and there are cofinally many such ordinals} $x$ {\rm (}$\mathcal{K}<a\in Ord \to \exists x\in Ord[a<x\land P_{I,n}(x)]${\rm ),}
 \item
 {\rm the axiom} `$\mathcal{K}$ {\rm is uncountable regular}'  {\rm is:}
\[
(\mathcal{K}>\omega)
\land
\forall\alpha<\mathcal{K}\forall f\in{}^{\alpha}\mathcal{K}\exists\beta<\mathcal{K}(f"\alpha\subset\beta)
\]
  {\rm and the axiom saying that} 
 $\forall B\subset \mathcal{K}[\mathcal{K}\in M(B)\to \exists\rho<\mathcal{K}(\rho\in M(B)\land Reg(\rho))]${\rm , which is codified
 by the following (\ref{eq:Kord}).}

\begin{equation}\label{eq:Kord}
\forall B\in L_{\mathcal{K}^{+}}
[B\subset \mathcal{K}
\to
\lnot\tau(B,\mathcal{K})
\to
\exists\rho<\mathcal{K}(\lnot\tau(B,\rho)\land Reg(\rho))
]
\end{equation}
{\rm where}
\begin{equation}\label{eq:thin}
\tau(B,\rho):\Leftrightarrow \exists C\subset \rho
[(C \mbox{{ \rm is club}})^{\rho} \land
(B\cap C=\emptyset)]
\end{equation}
{\rm and} $(C \mbox{{ \rm is club}})^{\rho}$ {\rm is a formula saying that} $C$ {\rm is a club subset of} $\rho$.

{\rm Namely}
$\tau(B,\rho)$ {\rm says that the set} $B$ {\rm is thin, i.e., non-stationary in} $\rho$.

{\rm Note that} $ (C \mbox{{ \rm is club}})^{\rho} \land (B\cap C=\emptyset)$ {\rm is stratified with respec to} $B,C${\rm , and}
$\tau(B,\rho)$ {\rm is stratified with respec to} $B$.
 \eenu

\eenu

\end{definition}

The following Lemma \ref{lem:regularset} is seen as in \cite{liftupZF}.

\begin{lemma}\label{lem:regularset}
$\mbox{{\rm T}}(\mathcal{K},I):=\bigcup_{n\in\omega}\mbox{{\rm T}}(\mathcal{K},I,n)$ is equivalent to the set theory 
${\sf ZF}+(V=L)+(\mathcal{K}\mbox{ {\rm is weakly compact)}}$.
\end{lemma}

\section{Operator controlled derivations for weakly compact cardinals}\label{sect:controlledOme}

In this section,
operator controlled derivations are first introduced,
and
inferences $(\mbox{{\bf Ref}}_{\mathcal{K}})$
for $\Pi^{1}_{1}$-indescribability are then eliminated from operator controlled derivations of $\Sigma^{1}_{2}$-sentences $\varphi^{V_{\mathcal{K}}}$ over $\mathcal{K}$.

In what follows $n$ denotes a fixed positive integer.
We tacitly assume that any ordinal is in $\mathcal{H}(n)$.

For $\alpha<^{\varepsilon}I=\langle 1,0\rangle$, $L_{\alpha}=\{\langle 0,x\rangle:x\in L_{(\alpha)_{1}}\}$.
$L_{I}=\{\langle 0,x\rangle :x\in L\}=\bigcup_{\alpha<^{\varepsilon}I}L_{\alpha}$
denotes the universe.
Both $(L_{I},\in^{\varepsilon})\models A$ and `$A$ is true' are synonymous with $A$.

\subsection{An intuitionistic fixed point theory $\mbox{FiX}^{i}(\mbox{{\sf ZFLK}}_{n})$}\label{subsec:intfixZFL}

For the fixed positive integer $n$,
$\mbox{{\sf ZFLK}}_{n}$ denotes the set theory 
${\sf ZF}+(V=L)+(\mathcal{K}\in Mh_{n}^{\omega_{n}(I+1)})$ in the language $\{\in,\mathcal{K}\}$
with an individual constant $\mathcal{K}$.
Let us also denote the set theory
${\sf ZF}+(V=L)+(\mathcal{K}\mbox{ is weakly compact})$ in the language $\{\in,\mathcal{K}\}$
by {\sf ZFLK}.

To analyze the theory {\sf ZFLK}, we need to handle
the relation
$(\mathcal{H}_{\gamma}[\Theta_{0}],\Theta,\kappa,n)\vdash^{a}_{b}\Gamma$ defined in subsection \ref{subsec:operatorcont},
where $n$ is the fixed integer, 
$\gamma,\kappa,a,b$ are codes of ordinals with $a<^{\varepsilon}\omega_{n}(I+1)$, $b<^{\varepsilon}I\oplus\omega$ and $\kappa\leq^{\varepsilon} I$ the code of a regular ordinal,
$\Theta_{0},\Theta$ are finite subsets of $L_{I}$ and $\Gamma$ a sequent, i.e., a finite set of sentences.
Usually the relation is defined by recursion on `ordinals' $a$, but such a recursion is not available in 
$\mbox{{\sf ZFLK}}_{n}$
since $a$ may be larger than $I$.
Instead of the recursion, the relation is defined for each $n<\omega$, as a fixed point,
\begin{equation}\label{eq:fixH}
H_{n}(\gamma,\Theta_{0},\Theta,\kappa,a,b,\Gamma)\Leftrightarrow (\mathcal{H}_{\gamma,n}[\Theta_{0}],\Theta,\kappa,n)\vdash^{a}_{b}\Gamma
\end{equation}
In this way the whole proof in this section is formalizable in 
an intuitionistic fixed point theory $\mbox{FiX}^{i}(\mbox{{\sf ZFLK}}_{n})$ over $\mbox{{\sf ZFLK}}_{n}$.

Throughout this section we work in an intuitionistic fixed point theory 
$\mbox{FiX}^{i}(\mbox{{\sf ZFLK}}_{n})$ over $\mbox{{\sf ZFLK}}_{n}$.
The intuitionistic theory $\mbox{FiX}^{i}(\mbox{{\sf ZFLK}}_{n})$ is introduced in \cite{liftupZF},
and shown to be a conservative extension of $\mbox{{\sf ZFLK}}_{n}$.
Let us reproduce definitions and results on $\mbox{FiX}^{i}(\mbox{{\sf ZFLK}}_{n})$ here.

Fix an $X$-strictly positive formula $\mathcal{Q}(X,x)$ in the language $\{\in,\mathcal{K},=,X\}$ with an extra unary predicate symbol $X$.
In $\mathcal{Q}(X,x)$ the predicate symbol $X$ occurs only strictly positive.
This means that the predicate symbol $X$ does not occur in the antecedent $\varphi$ of implications $\varphi\to\psi$ 
nor in the scope of negations $\lnot$ in $\mathcal{Q}(X,x)$.
The language of $\mbox{FiX}^{i}(\mbox{{\sf ZFLK}}_{n})$ is $\{\in,\mathcal{K},=,Q\}$ with a fresh unary predicate symbol $Q$.
The axioms in $\mbox{FiX}^{i}(\mbox{{\sf ZFLK}}_{n})$ consist of the following:
\begin{enumerate}
\item
All provable sentences in $\mbox{{\sf ZFLK}}_{n}$ (in the language $\{\in,\mathcal{K},=\}$).

\item
Induction schema for any formula $\varphi$ in $\{\in,\mathcal{K},=,Q\}$:
\begin{equation}\label{eq:Qind}
\forall x(\forall y\in x\,\varphi(y)\to\varphi(x))\to\forall x\,\varphi(x)
\end{equation}

\item
Fixed point axiom:
\[
\forall x[Q(x)\leftrightarrow \mathcal{Q}(Q,x)]
.\]
\end{enumerate}

The underlying logic in $\mbox{FiX}^{i}(\mbox{{\sf ZFLK}}_{n})$ is defined to be the intuitionistic (first-order predicate) logic (with equality).

(\ref{eq:Qind}) yields the following Lemma \ref{lem:vepsfix}.

\begin{lemma}\label{lem:vepsfix}
Let $<^{\varepsilon}$ denote a $\Delta_{1}$-predicate 
as described before Theorem \ref{th:resolutionpi3}.
For each $n<\omega$ and each formula $\varphi$ in $\{\in,\mathcal{K},=,Q\}$,
\[
\mbox{{\rm FiX}}^{i}(\mbox{{\sf ZFLK}}_{n})\vdash\forall x(\forall y<^{\varepsilon}x\,\varphi(y)\to\varphi(x)) \to 
\forall x<^{\varepsilon}\omega_{n}(I+1)\varphi(x)
.\]
\end{lemma}

The following Theorem \ref{th:consvintfix} is seen as in \cite{intfix, liftupZF}.

\begin{theorem}\label{th:consvintfix}
$\mbox{{\rm FiX}}^{i}(\mbox{{\sf ZFLK}}_{n})$ is a conservative extension of $\mbox{{\sf ZFLK}}_{n}$.
\end{theorem}

In what follows we work in $\mbox{FiX}^{i}(\mbox{{\sf ZFLK}}_{n})$ for a fixed integer $n$.

\subsection{Classes of sentences}

$\mathcal{K}\in L=L_{I}=\bigcup_{\alpha\in Ord}L_{\alpha}$ denotes a transitive and wellfounded model of ${\sf ZF}+(V=L)$, 
where
$L_{\alpha+1}$ is the set of $L_{\alpha}$-definable subsets of $L_{\alpha}$.
$Ord$ denotes the class of all ordinals in $L$, and $I$ the least ordinal not in $L$,
while $Ord^{\varepsilon}$ denotes the codes of ordinals less than $\omega_{n}(I+1)$.

\begin{definition}
{\rm For} $a\in L${\rm ,} $\mbox{{\rm rk}}_{L}(a)$ {\rm denotes the} $L$-rank {\rm of} $a$.
\[
\mbox{{\rm rk}}_{L}(a):=\min\{\alpha\in Ord: a\in L_{\alpha+1}\}
.\]
\end{definition}
If $a\in b\in L$, then $a\in b\subset L_{\beta}$ for $\beta=\mbox{{\rm rk}}_{L}(b)$ and $a\in L_{\beta}$.
Hence $\mbox{{\rm rk}}_{L}(a)<\beta=\mbox{{\rm rk}}_{L}(b)$.

The language $\mathcal{L}_{c}$ is obtained from the language $\{\in, P,P_{I,n},Reg,\mathcal{K}\}$ 
by adding names(individual constants) $c_{a}$
of each set $a\in L$.
$c_{a}$ is identified with $a$.

Then {\it formulae\/} in $\mathcal{L}_{c}$
is defined as usual.
Unbounded quantifiers $\exists x,\forall x$ are denoted by $\exists x\in L_{I},\forall x\in L_{I}$, resp.

For formulae $A$ in $\mathcal{L}_{c}$, 
${\sf qk}(A)$ denotes the finite set of $L$-ranks $\mbox{{\rm rk}}_{L}(a)$ 
of sets $a$ which are bounds of `bounded' quantifiers $\exists x\in a,\forall x\in a$
occurring in $A$.
Moreover
${\sf k}(A)$ denotes the set of $L$-ranks of sets occurring in $A$,
while ${\sf k}^{E}(A)$ denotes the set of $L$-ranks of sets occurring in an unstratifed position in $A$.
Both ${\sf k}(A)$ and ${\sf k}^{E}(A)$ are defined to include $L$-ranks of  bounds of `bounded' quantifiers.
Thus ${\sf qk}(A)\subset{\sf k}^{E}(A)\subset{\sf k}(A)\leq I$.
By definition we set $0\in{\sf qk}(A)$.

In the following definition, $Var$ denotes the set of variables and set $\mbox{{\rm rk}}_{L}(x):=0$ for variables $x\in Var$.

\begin{definition}
\benu

\item
${\sf k}(\lnot A)={\sf k}(A)$ {\rm and similarly for} ${\sf k}^{E},{\sf qk}$.
\item
${\sf qk}(M)=\{0\}$ {\rm for any literal} $M$.
\item
${\sf k}^{E}(M)={\sf k}(M)=\{\mbox{{\rm rk}}_{L}(t): t\in\vec{t}\,\}\cup\{0\}$
{\rm for literals} $Q(\vec{t}\,)$ {\rm with predicates} $Q\in\{P,P_{I,n},Reg\}$.
\item
${\sf k}(t\in s)=\{\mbox{{\rm rk}}_{L}(t),\mbox{{\rm rk}}_{L}(s), 0\}$ {\rm and}
${\sf k}^{E}(t\in s)=\{\mbox{{\rm rk}}_{L}(t), 0\}$.
\item
${\sf k}(A_{0}\lor A_{1})={\sf k}(A_{0})\cup{\sf k}(A_{1})$ {\rm and similarly for} ${\sf k}^{E},{\sf qk}$.
\item
{\rm For} $t\in L_{I}\cup\{L_{I}\}\cup Var$,
${\sf k}(\exists x\in t\, A(x))=\{\mbox{{\rm rk}}_{L}(t)\}\cup{\sf k}(A(x))$ {\rm and similarly for} ${\sf k}^{E},{\sf qk}$.
\eenu

\end{definition}

For example
${\sf k}^{E}(a\in b)=\{\mbox{{\rm rk}}_{L}(a),0\}$, and
${\sf qk}(\exists x\in a\, A(x))=\{\mbox{{\rm rk}}_{L}(a)\}\cup{\sf qk}(A(x))$.

\begin{definition}\label{df:fmlclasses}
\benu

\item
$A\in\Delta_{0}$ {\rm iff there exists a} $\Delta_{0}${\rm -formula} $\theta[\vec{x}\,]$ 
{\rm in the language} $\{\in\}$ {\rm and terms}
 $\vec{t}$ {\rm such that} $A\equiv\theta[\vec{t}\,]$.
{\rm This means that} $A$ {\rm is bounded, and the predicates} $P,P_{I,n},Reg$
{\rm do not occur in} $A$.

\item
{\rm Putting} $\Sigma_{0}:=\Pi_{0}:=\Delta_{0}${\rm , the classes}
$\Sigma_{m}$ {\rm and} $\Pi_{m}$ {\rm of formulae in the language} $\{\in\}$ {\rm with terms}
{\rm are defined as usual using quantifiers}
$\exists x\in L_{I},\forall x\in L_{I}${\rm , where by definition}
$\Sigma_{m}\cup\Pi_{m}\subset\Sigma_{m+1}\cap\Pi_{m+1}$.

{\rm Each formula in} $\Sigma_{m}\cup\Pi_{m}$ {\rm is in prenex normal form with alternating unbounded quantifiers and}
$\Delta_{0}${\rm -matrix.}

\item
$A\in\Delta_{0}(\lambda)$ {\rm iff there exists a} $\Delta_{0}${\rm -formula} $\theta[\vec{x}\,]$ 
{\rm in the language} $\{\in\}$ {\rm and terms}
 $\vec{t}$ {\rm such that} $A\equiv\theta[\vec{t}\,]$ {\rm and}
 ${\sf k}(A)<\lambda$.

\item 
$A\in \Sigma_{1}(\lambda)$ {\rm iff either} $A\in\Delta_{0}(\lambda)$ {\rm or}
$A\equiv\exists x\in L_{\lambda}\, B$ {\rm with} $B\in\Delta_{0}(\lambda)$.

{\rm Note that} $\Sigma(\lambda)\subset\Delta_{0}$ {\rm for any} $\lambda<I$.

\item
{\rm The class of sentences} $\Sigma_{m}(\lambda),\Pi_{m}(\lambda)\, (m<\omega)$ {\rm are defined as usual.}

\item
$\Sigma^{1}_{0}(\lambda)$ {\rm denotes the set of first-order formulae on} $L_{\lambda}${\rm , i.e.,}
$\Sigma^{1}_{0}(\lambda):=\bigcup_{m\in\omega}\Sigma_{m}(\lambda)$.
\eenu

\end{definition}

Note that the predicates $P,P_{I,n},Reg$ do not occur in $\Sigma_{m}$-formulae nor in $\Sigma^{1}_{0}(\lambda)$-formulae.

\begin{definition}
{\rm A set} $\Sigma^{\Sigma_{n+1}}(\lambda)$ {\rm of sentences is defined recursively as follows.}
 \benu

 \item
$\Sigma_{n+1}\subset\Sigma^{\Sigma_{n+1}}(\lambda)$.

 \item
 {\rm Each literal including} $Reg(a), P(a,b,c),P_{I,n}(a)$ {\rm and their negations 
 is in} $\Sigma^{\Sigma_{n+1}}(\lambda)$.

 \item
 $\Sigma^{\Sigma_{n+1}}(\lambda)$ {\rm is closed under propositional connectives} $\lor,\land$.

 \item
 {\rm Suppose}
 $\forall x\in b\, A(x)\not\in\Delta_{0}${\rm . Then}
 $\forall x\in b\, A(x)\in \Sigma^{\Sigma_{n+1}}(\lambda)$ {\rm iff} $A(\emptyset)\in \Sigma^{\Sigma_{n+1}}(\lambda)$ {\rm and}
 $\mbox{{\rm rk}}_{L}(b)<\lambda$.

 \item
 {\rm Suppose}
 $\exists x\in b\, A(x)\not\in\Delta_{0}${\rm . Then}
 $\exists x\in b\, A(x)\in \Sigma^{\Sigma_{n+1}}(\lambda)$ {\rm iff} $A(\emptyset)\in \Sigma^{\Sigma_{n+1}}(\lambda)$ {\rm and}
 $\mbox{{\rm rk}}_{L}(b)\leq\lambda$.
 
 \eenu

 \end{definition}

\begin{definition}\label{df:domFfml}
{\rm Let us extend the domain} $dom(F^{\Sigma_{1}}_{x\cup\{\kappa\}})=\mbox{{\rm Hull}}_{\Sigma_{1}}^{I}(x\cup\{\kappa\})$
{\rm of Mostowski collapse to formulae.}
\[
dom(F^{\Sigma_{1}}_{x\cup\{\kappa\}})=\{A\in\Sigma_{1}\cup\Pi_{1}: {\sf k}(A)\subset\mbox{{\rm Hull}}_{\Sigma_{1}}^{I}(x\cup\{\kappa\})\}
.\]
{\rm For} $A\in dom(F^{\Sigma_{1}}_{x\cup\{\kappa\}})$,
$F^{\Sigma_{1}}_{x\cup\{\kappa\}}" A$ {\rm denotes the result of replacing each constant} $\gamma$ {\rm by} 
$F^{\Sigma_{1}}_{x\cup\{\kappa\}}(\gamma)${\rm , 
 each unbounded existential quantifier} $\exists z\in L_{I}$ {\rm by} $\exists z\in L_{F^{\Sigma_{1}}_{x\cup\{\kappa\}}(I)}${\rm ,
and each unbounded universal quantifier} $\forall z\in L_{I}$ {\rm by} $\forall z\in L_{F^{\Sigma_{1}}_{x\cup\{\kappa\}}(I)}$.

{\rm For sequent, i.e., finite set of sentences} $\Gamma\subset dom(F^{\Sigma_{1}}_{x\cup\{\kappa\}})${\rm , put}
 $F^{\Sigma_{1}}_{x\cup\{\kappa\}}"\Gamma=\{F^{\Sigma_{1}}_{x\cup\{\kappa\}}" A: A\in\Gamma\}$.
 
 {\rm Likewise the domain} $dom(F^{\Sigma_{n}}_{x})=\mbox{{\rm Hull}}_{\Sigma_{n}}^{I}(x)$
{\rm is extended to}
\[
dom(F^{\Sigma_{n}}_{x})=\{A\in\Sigma_{n}\cup\Pi_{n}: {\sf k}(A)\subset\mbox{{\rm Hull}}_{\Sigma_{n}}^{I}(x)\}
\]
{\rm and for formula} $A\in dom(F^{\Sigma_{n}}_{x})$,
$F^{\Sigma_{n}}_{x}" A$, {\rm and sequent} $\Gamma\subset dom(F^{\Sigma_{n}}_{x})$,
 $F^{\Sigma_{n}}_{x}"\Gamma$
 {\rm are defined similarly.}
 
\end{definition}

\begin{proposition}\label{prp:domFfml}
For $F=F^{\Sigma_{1}}_{x\cup\{\kappa\}}, F^{\Sigma_{n}}_{x}$ and $A\in dom(F)$
\[
L_{I}\models A\leftrightarrow F" A
.\]
\end{proposition}

The assignment of disjunctions and conjunctions to sentences is defined as in
\cite{liftupZF}.

\begin{definition}\label{df:assigndc}
\benu

\item\label{df:assigndc0}
{\rm If} $M$ {\rm is one of the literals} $a\in b,a\not\in b${\rm , then for} $J:=0$
\[
M:\simeq
\left\{
\begin{array}{ll}
\bigvee(A_{\iota})_{\iota\in J} & \mbox{{\rm if }} M \mbox{ {\rm is false (in }} L_{I}\mbox{{\rm )}}
\\
\bigwedge(A_{\iota})_{\iota\in J} &  \mbox{{\rm if }} M \mbox{ {\rm is true}}
\end{array}
\right.
\]

\item
$(A_{0}\lor A_{1}):\simeq\bigvee(A_{\iota})_{\iota\in J}$
{\rm and}
$(A_{0}\land A_{1}):\simeq\bigwedge(A_{\iota})_{\iota\in J}$
{\rm for} $J:=2$.

\item

\[
Reg(a) :\simeq
\bigvee(a=a)_{\iota\in J} 
\mbox{ {\rm and }}
\lnot Reg(a):\simeq
\bigwedge(a\neq a)_{\iota\in 1} 
\]
{\rm with}
\[
J:=
\left\{
\begin{array}{ll}
1& \mbox{{\rm if }} a\in R
\\
0 & \mbox{{\rm otherwise}}
\end{array}
\right.
.\]

\item
\[
P(a,b,c) :\simeq
\bigvee(a=a)_{\iota\in J} 
\mbox{ {\rm and }}
\lnot P(a,b,c) :\simeq
\bigwedge(a\neq a)_{\iota\in J} 
\]

{\rm with}
\[
J:=
\left\{
\begin{array}{ll}
1 & \mbox{{\rm if }}  
a\in R \,\&\, \exists \alpha\in Ord_{\varepsilon}[b=\Psi_{a,n}\alpha\,\&\,
\alpha\in\mathcal{H}_{\alpha}(b)\,\&\,  c=F^{\Sigma_{1}}_{b\cup\{a\}}(I)]
\\
0 & \mbox{{\rm otherwise}}
\end{array}
\right.
.\]

\item
\[
P_{I,n}(a) :\simeq
\bigvee(a=a)_{\iota\in J} 
\mbox{ {\rm and }}
\lnot P_{I,n}(a) :\simeq
\bigwedge(a\neq a)_{\iota\in J} 
\]

{\rm with}
\[
J:=
\left\{
\begin{array}{ll}
1 & \mbox{{\rm if }}  
 \exists \alpha\in Ord_{\varepsilon}[a=\Psi_{I,n}\alpha\,\&\,
\alpha\in\mathcal{H}_{\alpha}(a)]
\\
0 & \mbox{{\rm otherwise}}
\end{array}
\right.
.\]

\item
{\rm Let} $(\exists z\in b \, \theta[z])\in\Sigma_{n}$ {\rm for} $b\in L_{I}\cup\{L_{I}\}${\rm , and}
$(\exists z\in b \, \theta[z])\not\in\Sigma^{1}_{0}(\mathcal{K}^{+})$.
{\rm Then for the set}
\begin{equation}\label{eq:dfmu}
\mu z\in b \, \theta[z] :=
\min_{<_{L}}\{d : (d\in b \land \theta[d]) \lor (\lnot\exists z\in b\, \theta[z]\land d=0)\}
\end{equation}
{\rm with a canonical well ordering} $<_{L}$ {\rm on} $L$
{\rm , and}
$J=\{d\}$
\begin{eqnarray}
\exists z\in b\, \theta[z] & :\simeq & \bigvee(d\in b\land \theta[d])_{d\in J}
\label{eq:sigpimu}
\\
\forall z\in b\, \lnot\theta[z] & :\simeq & \bigwedge(d\in b \to \lnot\theta[d)_{d\in J}
\nonumber
\end{eqnarray}
{\rm where} $d\in b$ {\rm denotes a true literal, e.g.,} $d\not\in d$ {\rm when} $b=L_{I}$.

{\rm This case is applied only when} $\exists z\in b\, \theta[z]$ {\rm is a formula in} $\{\in\}\cup L_{I}${\rm , and}
$(\exists z\in b \, \theta[z])\in\Sigma_{n}$ {\rm but}
$(\exists z\in b \, \theta[z])\not\in\Sigma^{1}_{0}(\mathcal{K}^{+})$.

\item
{\rm Otherwise set for} $a\in L_{I}\cup\{L_{I}\}$
\[
\exists x\in a\, A(x):\simeq\bigvee(A(b))_{b\in J}
\mbox{ {\rm and }}
\forall x\in a\, A(x):\simeq\bigwedge(A(b))_{b\in J}
\]
{\rm for} 
\[
J:=\{b: b\in a\}
.\]
{\rm This case is applied if one of the predicates} $P,P_{I,n},Reg$ {\rm occurs in} $\exists x\in a\, A(x)$
{\rm , or} $(\exists x\in a\, A(x))\not\in\Sigma_{n}${\rm , or} $(\exists x\in a\, A(x))\in\Sigma^{1}_{0}(\mathcal{K}^{+})$.
\eenu

\end{definition}

In particular we have
\begin{eqnarray*}
\lnot\tau(B,\mathcal{K})
& :\simeq &
\bigwedge\{
(C\not\subset\mathcal{K})\lor
\lnot (C\mbox{{ \rm is club}})^{\mathcal{K}}\lor
(B\cap C\neq\emptyset)
:
C\in L_{\mathcal{K}^{+}}\}
\\
\tau(B,\mathcal{K})
& :\simeq &
\bigvee\{
(C\subset\mathcal{K})\land
 (C\mbox{{ \rm is club}})^{\mathcal{K}}\land
(B\cap C=\emptyset)
:
C\in L_{\mathcal{K}^{+}}\}
\end{eqnarray*}
{\rm where}
\begin{equation}
\renewcommand{\theequation}{\ref{eq:thin}}
\tau(B,\rho):\Leftrightarrow \exists C\subset\rho(
[(C \mbox{{ \rm is club}})^{\rho} \land
(B\cap C=\emptyset)]
\end{equation}
\addtocounter{equation}{-1}

The definition of the rank $\mbox{{\rm rk}}(A)$ of sentences $A$ in \cite{liftupZF}
is slightly changed as follows.
The rank $\mbox{{\rm rk}}(A)$ of sentences $A$ is defined by recursion on the  number of symbols occurring in 
$A$.

\begin{definition}\label{df:rank}
\benu

\item\label{df:rank1}
$\mbox{{\rm rk}}(\lnot A):=\mbox{{\rm rk}}(A)$.

\item\label{df:rank2}
$\mbox{{\rm rk}}(a\in b):=\mbox{{\rm rk}}(a\not\in b):=0$.

\item\label{df:rank5-1}
$\mbox{{\rm rk}}(Reg(\alpha)):=\mbox{{\rm rk}}(P(\alpha,\beta,\gamma)):=\mbox{{\rm rk}}(P_{I,n}(\alpha)):=1$.

\item\label{df:rank6}
$\mbox{{\rm rk}}(A_{0}\lor A_{1}):=\max\{\mbox{{\rm rk}}(A_{0}),\mbox{{\rm rk}}(A_{1})\}+1$.

\item\label{df:rank7}
$
\mbox{{\rm rk}}(\exists x\in a\, A(x)):=
\max\{\omega\alpha, \mbox{{\rm rk}}(A(\emptyset))+2\}
$
{\rm for} $\alpha=\mbox{{\rm rk}}_{L}(a)$.
\eenu

\end{definition}

\begin{proposition}\label{lem:rank}
Let $A\simeq\bigvee(A_{\iota})_{\iota\in J}$ or $A\simeq\bigwedge(A_{\iota})_{\iota\in J}$.
\benu

\item\label{lem:rank15}
$A\in \Sigma^{\Sigma_{n+1}}(\lambda)\Rightarrow\forall\iota\in J(A_{\iota}\in \Sigma^{\Sigma_{n+1}}(\lambda))$.

\item\label{prp:rksig3}
For an ordinal $\lambda\leq I$ with $\omega\lambda=\lambda$,
$
\mbox{{\rm rk}}(A)<\lambda \Rightarrow A\in\Sigma^{\Sigma_{n+1}}(\lambda)
$.

\item\label{lem:rank0}
$\mbox{{\rm rk}}(A)<I+\omega$.

\item\label{lem:rank1}
$\mbox{{\rm rk}}(A)$ is in the Skolem hull of $\omega{\sf qk}(A)\cup\{0,1\}$ under the addition
with $\omega{\sf qk}(A)=\{\omega\alpha:\alpha\in{\sf qk}(A)\}$.

\item\label{lem:rank2}
$\forall\iota\in J(\mbox{{\rm rk}}(A_{\iota})<\mbox{{\rm rk}}(A))$.

\eenu

\end{proposition}
\bprf

\ref{lem:rank}.\ref{lem:rank2}.
This is seen from the fact that
$a\in b\in L \Rightarrow \mbox{{\rm rk}}_{L}(a)<\mbox{{\rm rk}}_{L}(b)$.
\eprf

\subsection{Operator controlled derivations}\label{subsec:operatorcont}

$\kappa,\lambda,\sigma,\pi$ ranges over $R^{+}$.

Let $\mathcal{H}$ be an operator, $\Theta$ a finite set of ordinals, 
$\kappa \in R^{+}$, $\Gamma$ a sequent, $a\in Ord^{\varepsilon}$ and 
$b<I+\omega$.
We define a relation $(\mathcal{H},\Theta,\kappa,n)\vdash^{a}_{b}\Gamma$, which is read `there exists an infinitary derivation
of $\Gamma$ which is $(\kappa,n)${\it -controlled\/} by $\mathcal{H}$ and $\Theta$, and 
 whose height is at most $a$ and its cut rank is less than $b$'.

Recall that $R$ denotes the set of uncountable cardinals $\rho$ such that $\mathcal{K}<\rho<I$, and
$\lambda>\mathcal{K}$
in the inference rules $(\mbox{{\bf P}}_{\lambda})$ and $(\mbox{{\bf F}}^{\Sigma_{1}}_{x\cup\{\lambda\}})$.

Sequents are finite sets of sentences, and inference rules are formulated in one-sided sequent calculus.

\begin{definition}\label{df:sfkEK}
\[
{\sf k}^{E}_{\mathcal{K}}(A):=
\left\{
\begin{array}{ll}
{\sf k}^{E}(A) & \mbox{{\rm if }} A\in\Sigma^{1}_{0}(\mathcal{K}^{+})
\\
{\sf k}(A) &  \mbox{{\rm otherwise}}
\end{array}
\right.
\]
\end{definition}

\begin{definition}\label{df:controlderreg}
$(\mathcal{H},\Theta,\kappa,n)\vdash^{a}_{b}\Gamma$ {\rm holds if} 
\begin{equation}\label{eq:controlder}
{\sf k}^{E}_{\mathcal{K}}(\Gamma):=\bigcup\{{\sf k}^{E}_{\mathcal{K}}(A):A\in\Gamma\}\subset\mathcal{H}:=\mathcal{H}(\emptyset)
\,\&\,
a\in\mathcal{H}[\Theta]
\end{equation}
{\rm and one of the following
cases holds:}

\benu

\item
 $A\simeq\bigvee\{A_{\iota}: \iota\in J\}$, $A\in\Gamma$ {\rm and for an} $\iota\in J$, $a(\iota)<a$ {\rm and}
$\mbox{{\rm rk}}_{L}(\iota)<\kappa 
\Rightarrow \mbox{{\rm rk}}_{L}(\iota)< a
$

\[
\infer[(\bigvee)]{(\mathcal{H},\Theta,\kappa,n)\vdash^{a}_{b}\Gamma}
{(\mathcal{H},\Theta,\kappa,n)\vdash^{a(\iota)}_{b}\Gamma,A_{\iota}}
\]

\item
$A\simeq\bigwedge\{A_{\iota}: \iota\in J\}$, $A\in\Gamma$ {\rm and} $a(\iota)<a$ {\rm for any} $\iota\in J$
\[
\infer[(\bigwedge)]{(\mathcal{H},\Theta,\kappa,n)\vdash^{a}_{b}\Gamma}
{\{(\mathcal{H}[\{\mbox{{\rm rk}}_{L}(\iota)\}],\Theta,\kappa,n)\vdash^{a(\iota)}_{b}\Gamma,A_{\iota}:\iota\in J\}}
\]

\item

$\mbox{{\rm rk}}(C)<b$ {\rm and an} $a_{0}<a$
\[
\infer[(cut)]{(\mathcal{H},\Theta,\kappa,n)\vdash^{a}_{b}\Gamma}
{(\mathcal{H},\Theta,\kappa,n)\vdash^{a_{0}}_{b}\Gamma,\lnot C & (\mathcal{H},\Theta,\kappa,n)\vdash^{a_{0}}_{b}C,\Gamma}
\]

\item
$\alpha<\lambda\in R$ {\rm and} $\{\exists x<\lambda\exists y<\lambda[\alpha<x \land P(\lambda,x,y)]\}\cup\Gamma_{0}=\Gamma$
\[
\infer[(\mbox{{\bf  P}}_{\lambda})]{\exists x<\lambda\exists y<\lambda[\alpha<x \land P(\lambda,x,y)],\Gamma_{0}}{}
\]

\item

{\rm Let} $\lambda\in R$
{\rm and} $x\in\mathcal{H}[\Theta]$ {\rm where for some} $b$
\[
 x=
\Psi_{\lambda,n}b
.\]

{\rm If} 
$\Gamma=\Lambda\cup (F^{\Sigma_{1}}_{x\cup\{\lambda\}}"\Gamma_{0})$, 
$\Gamma_{0}\subset\Sigma_{1}$, $a_{0}<a$
{\rm and}
\[
{\sf k}(\Gamma_{0})\subset\mbox{{\rm Hull}}_{\Sigma_{1}}^{I}((\mathcal{H}\cap x)\cup\{\lambda\})
\]
{\rm  then}
\[
\infer[(\mbox{{\bf F}}^{\Sigma_{1}}_{x\cup\{\lambda\}})]{(\mathcal{H},\Theta,\kappa,n)\vdash^{a}_{b}
\Lambda,F^{\Sigma_{1}}_{x\cup\{\lambda\}}"\Gamma_{0}}{(\mathcal{H},\Theta,\kappa,n)\vdash^{a_{0}}_{b}\Lambda,\Gamma_{0}}
\]
{\rm where} $F^{\Sigma_{1}}_{x\cup\{\lambda\}}$ {\rm denotes the Mostowski collapse}
$F^{\Sigma_{1}}_{x\cup\{\lambda\}}:
 \mbox{{\rm Hull}}^{I}_{\Sigma_{1}}(x\cup\{\lambda\})\leftrightarrow L_{F^{\Sigma_{1}}_{x\cup\{\lambda\}}(I)}$.

\item
$\alpha<I$ {\rm and} $\{\exists x<I[\alpha<x \land P_{I,n}(x)]\}\cup\Gamma_{0}=\Gamma$
\[
\infer[(\mbox{{\bf  P}}_{I,n})]{ \exists x<I[\alpha<x \land P_{I,n}(x)],\Gamma_{0}}{}
\]

\item

{\rm Let}
\[
x=\Psi_{I,n}b\in\mathcal{H}[\Theta]
.\]
{\rm If} 
$\Gamma=\Lambda\cup (F^{\Sigma_{n}}_{x}"\Gamma_{0})$, 
$\Gamma_{0}\subset\Sigma_{n}$, $a_{0}<a$
{\rm and}
\[
{\sf k}(\Gamma_{0})\subset\mbox{{\rm Hull}}_{\Sigma_{n}}^{I}(\mathcal{H}\cap x)
\]
{\rm  then}
\[
\infer[(\mbox{{\bf F}}^{\Sigma_{n}}_{x})]{(\mathcal{H},\Theta,\kappa,n)\vdash^{a}_{b}\Lambda,F^{\Sigma_{n}}_{x}"\Gamma_{0}}
{(\mathcal{H},\Theta,\kappa,n)\vdash^{a_{0}}_{b}\Lambda,\Gamma_{0}}
\]
{\rm where} $F^{\Sigma_{n}}_{x}$ {\rm denotes the Mostowski collapse}
$F^{\Sigma_{n}}_{x}: \mbox{{\rm Hull}}^{I}_{\Sigma_{n}}(x)\leftrightarrow L_{F^{\Sigma_{n}}_{x}(I)}$.

\item

{\rm If} 
$\max\{a_{\ell},a_{r}\}<a$, 
 {\rm and}
$B\subset \mathcal{K}$,
$B\in\mbox{{\rm Hull}}^{I}_{\Sigma_{1}}(\{\mathcal{K},\mathcal{K}^{+}\})$,
{\rm then}

\[
\infer[(\mbox{{\bf Ref}}_{\mathcal{K}})]{(\mathcal{H},\Theta,\kappa,n)\vdash^{a}_{b}\Gamma}
{
(\mathcal{H},\Theta,\kappa,n)\vdash^{a_{\ell}}_{b}\Gamma, \lnot\tau(B,\mathcal{K})
&
(\mathcal{H},\Theta,\kappa,n)\vdash^{a_{r}}_{b}\Gamma, 
\forall\rho<\mathcal{K}\,\tau(B,\rho)
}
\]

{\rm where}
\begin{equation}
\renewcommand{\theequation}{\ref{eq:thin}}
\tau(B,\rho):\Leftrightarrow \exists C\subset\rho
[(C \mbox{{ \rm is club}})^{\rho} \land
(B\cap C=\emptyset)]
\end{equation}
\addtocounter{equation}{-1}
{\rm which is stratified with respec to} $B$.

\eenu

\end{definition}

An inspection to Definition \ref{df:controlderreg}
shows that there exists a strictly positive formula $H_{n}$
such that the relation $(\mathcal{H}_{\gamma,n}[\Theta_{0}],\Theta,\kappa,n)\vdash^{a}_{b}\Gamma$
is a fixed point of $H_{n}$ as in (\ref{eq:fixH}).

In what follows the relation should be understood as a fixed point of $H_{n}$,
and recall that we are working in 
the intuitionistic fixed point theory $\mbox{FiX}^{i}({\sf ZFLK}_{n})$
over ${\sf ZFLK}_{n}$ defined in subsection \ref{subsec:intfixZFL}.

\begin{proposition}\label{prp:kapcontrolled}
$(\mathcal{H},\Theta,\kappa,n)\vdash^{a}_{b}\Gamma \,\&\, \lambda\leq\kappa \Rightarrow (\mathcal{H},\Theta,\lambda,n)\vdash^{a}_{b}\Gamma$.
\end{proposition}

We will state some lemmata for the operator controlled derivations with sketches of their proofs
since
these can be shown as in \cite{Buchholz} and \cite{liftupZF}.

In what follows by an operator we mean an $\mathcal{H}_{\gamma}[\Theta]$ for a finite set $\Theta$ of ordinals.

\[
(\mathcal{H},\kappa,n)\vdash^{a}_{b}\Gamma
 :\Leftrightarrow 
(\mathcal{H},\emptyset,\kappa,n)\vdash^{a}_{b}\Gamma
\]

\begin{lemma}\label{lem:tautology}{\rm (Tautology)}
\[
(\mathcal{H}[{\sf k}^{E}_{\mathcal{K}}(A)],I,n)\vdash^{I+2\footnotesize{\mbox{{\rm rk}}}(A)}_{0}\Gamma,\lnot A, A
.\]

\end{lemma}

\begin{lemma}\label{lem:delta0complete}{\rm ($\Delta_{0}(I)$-completeness)}
If $\Gamma\subset\Delta_{0}(I)$ and $\bigvee\Gamma$ is true, then
\[
(\mathcal{H}[{\sf k}^{E}_{\mathcal{K}}(\Gamma)],I,n)\vdash^{I+2\footnotesize{\mbox{{\rm rk}}}(\Gamma)}_{0}\Gamma
\]
where $\mbox{{\rm rk}}(\Gamma)=\mbox{{\rm rk}}(A_{0})\#\cdots\#\mbox{{\rm rk}}(A_{n})$ for $\Gamma=\{A_{0},\ldots,A_{n}\}$.
\end{lemma}

\begin{lemma}\label{lem:falsedelta0elim}{\rm (Elimination of false sentences)}
\\

Let $A$ be a false sentence, i.e., $L_{I}\not\models A$, such that 
${\sf k}(A)\subset\mbox{{\rm Hull}}_{\Sigma_{1}}^{I}((\mathcal{K}+1)\cup\{\mathcal{K}^{+}\})\cap\mathcal{K}^{+}$.
Then
\[
(\mathcal{H},\Theta,\kappa,n)\vdash^{a}_{b}\Gamma,A \Rightarrow
(\mathcal{H},\Theta,\kappa,n)\vdash^{a}_{b}\Gamma
.\]
\end{lemma}
\bprf

Consider the case when $A$ is a main formula of an $({\bf F}^{\Sigma_{1}}_{x\cup\{\mathcal{K}^{+}\}})$ with $x>\mathcal{K}$.
We have $F^{\Sigma_{1}}_{x\cup\{\mathcal{K}^{+}\}}(a)=a$ for any $a$ with $\mbox{{\rm rk}}_{L}(a)<x$.

We claim $F^{\Sigma_{1}}_{x\cup\{\mathcal{K}^{+}\}}"A\equiv A$.
Let $b\in{\sf k}(A)$. 
Then $\mbox{{\rm rk}}_{L}(b)\in\mbox{{\rm Hull}}_{\Sigma_{1}}^{I}((\mathcal{K}+1)\cup\{\mathcal{K}^{+}\})\cap\mathcal{K}^{+}\subset\mbox{{\rm Hull}}^{I}_{\Sigma_{1}}(x\cup\{\mathcal{K}^{+}\})\cap\mathcal{K}^{+}\subset x$.
Hence $F^{\Sigma_{1}}_{x\cup\{\mathcal{K}^{+}\}}(b)=b$.
\eprf

\begin{lemma}\label{th:embedreg}{\rm (Embedding)}\\
 For each axiom $A$ in $\mbox{{\rm T}}(\mathcal{K},I,n)$, there is an $m<\omega$ such that
 for any operator $\mathcal{H}$ 
 \[
 (\mathcal{H}[\{\mathcal{K}\}],I,n)\vdash_{I}^{I\cdot m}{\rm `} \mathcal{K} \mbox{ {\rm is uncountable regular}' } \to A
  .\]
\end{lemma}
\bprf

The axiom for $\Pi^{1}_{1}$-indescribability
\begin{equation}
\renewcommand{\theequation}{\ref{eq:Kord}}
\forall B\in L_{\mathcal{K}^{+}}
[B\subset \mathcal{K}
\to
\lnot\tau(B,\mathcal{K})
\to
\exists\rho<\mathcal{K}(\lnot\tau(B,\rho)\land Reg(\rho))
]
\end{equation}
\addtocounter{equation}{-1}
follows from the inference rule $(\mbox{{\bf Ref}}_{\mathcal{K}})$ and 
$\mbox{(\ref{eq:Kord})}
\simeq
\bigwedge(B\subset \mathcal{K} \to \lnot\tau(B,\mathcal{K})
\to
\exists\rho<\mathcal{K}(\lnot\tau(B,\rho)\land Reg(\rho))_{B\in L_{\mathcal{K}^{+}}}$
for
$B:=\mu B\in L_{\mathcal{K}^{+}}(B\subset \mathcal{K} \land \lnot\tau(B,\mathcal{K})
\land
\forall\rho<\mathcal{K}(Reg(\rho)\to \tau(B,\rho)))\in\mbox{{\rm Hull}}^{I}_{\Sigma_{1}}(\{\mathcal{K},\mathcal{K}^{+}\})$.
\eprf

\begin{lemma}\label{lem:inversionreg}{\rm (Inversion)}
\\
Let  $d=\mu z\in b\,A[\vec{c},z] $ for $(\exists z\in b\, A)\in \Sigma_{n}\setminus\Sigma^{1}_{0}(\mathcal{K}^{+})$.

\[
(\mathcal{H},\Theta,\kappa,n)\vdash^{a}_{b}\Gamma, \exists z\in b\,A[\vec{c},z] \Rightarrow
(\mathcal{H},\Theta,\kappa,n)\vdash^{a}_{b}\Gamma, d\in b\land A[\vec{c},d]
\]
and
\[
(\mathcal{H},\Theta,\kappa,n)\vdash^{a}_{b}\Gamma, \forall z\in b\,\lnot A[\vec{c},z] \Rightarrow
(\mathcal{H},\Theta,\kappa,n)\vdash^{a}_{b}\Gamma, d\in b \to \lnot A[\vec{c},d]
\]
\end{lemma}

\begin{lemma}\label{lem:reduction}{\rm (Reduction)}\\
Let  $C\simeq\bigvee(C_{\iota})_{\iota\in J}$.
\benu

\item\label{lem:reduction1}
Suppose
$C\not\in\{\exists x<\lambda\exists y<\lambda[\alpha<x \land P(\lambda,x,y)]: \alpha<\lambda\in R\}\cup\{\exists x<I[\alpha<x \land P_{I,n}(x)]:\alpha<I\}$.

Then
\[
(\mathcal{H},\Theta,\kappa,n)\vdash^{a}_{c}\Delta,\lnot C \,\&\, (\mathcal{H},\kappa,n)\vdash^{b}_{c}C,\Gamma \,\&\, 
\mathcal{K}\leq\mbox{{\rm rk}}(C)\leq c
\Rightarrow
(\mathcal{H},\Theta,\kappa,n)\vdash^{a+b}_{c}\Delta,\Gamma
\]
\item\label{lem:reduction2}
Assume $C\equiv(\exists x<\lambda\exists y<\lambda[\alpha<x \land P(\lambda,x,y)])$ for an $\alpha<\lambda\in R$ and $\beta\in\mathcal{H}_{\beta}$.

Then
\[
(\mathcal{H}_{\beta},\kappa,n)\vdash^{a}_{b}\Gamma,\lnot C 
\Rightarrow
(\mathcal{H}_{\beta+1},\kappa,n)\vdash^{a}_{b}\Gamma
\]

\item\label{lem:reduction3}
Assume $C\equiv(\exists x<I[\alpha<x \land P_{I,n}(x)])$ for an $\alpha<I$ and $\beta\in\mathcal{H}_{\beta}$.

Then
\[
(\mathcal{H}_{\beta},\kappa,n)\vdash^{a}_{b}\Gamma,\lnot C 
\Rightarrow
(\mathcal{H}_{\beta+1},\kappa,n)\vdash^{a}_{b}\Gamma
\]
\eenu

\end{lemma}

\begin{lemma}\label{lem:predcereg}{\rm (Predicative Cut-elimination)}
\benu

\item\label{lem:predcereg2}
$(\mathcal{H},\kappa,n)\vdash^{b}_{c+\omega^{a}}\Gamma
\,\&\, 
[c,c+\omega^{a}[\cap (\{\lambda+1:\lambda\in R\}\cup\{I\})=\emptyset
\,\&\, a\in\mathcal{H}
\Rightarrow (\mathcal{H},\kappa,n)\vdash^{\varphi ab}_{c}\Gamma$.
\item\label{lem:predcereg4}
For $\lambda\in R$,
$(\mathcal{H}_{\gamma},\kappa,n)\vdash^{b}_{\lambda+2}\Gamma \,\&\, \gamma\in\mathcal{H}_{\gamma}\,\&\,
 \Rightarrow 
(\mathcal{H}_{\gamma+b},\kappa,n)\vdash^{\omega^{b}}_{\lambda+1}\Gamma$.
\item\label{lem:predcereg5}
$(\mathcal{H}_{\gamma},\kappa,n)\vdash^{b}_{I+1}\Gamma \,\&\, \gamma\in\mathcal{H}_{\gamma}\,\&\,
 \Rightarrow 
(\mathcal{H}_{\gamma+b},\kappa,n)\vdash^{\omega^{b}}_{I}\Gamma$.
\item\label{lem:predcereg6}
$(\mathcal{H}_{\gamma},\kappa,n)\vdash^{b}_{c+\omega^{a}}\Gamma
\,\&\, 
\max\{a,b,c\}<I
\,\&\, a\in\mathcal{H}_{\gamma}
\Rightarrow (\mathcal{H}_{\gamma+\varphi ab},\kappa,n)\vdash^{\varphi ab}_{c}\Gamma$.
\eenu

\end{lemma}

\begin{definition}
{\rm For a formula} $\exists x\in d\, A$ {\rm and ordinals} $\lambda=\mbox{{\rm rk}}_{L}(d)\in R^{+}, \alpha$,
$(\exists x\in d\, A)^{(\exists\lambda\!\upharpoonright\!\alpha)}$ {\rm denotes the result of restricting the} {\it outermost existential quantifier\/}
$\exists x\in d$ {\rm to} $\exists x\in L_{\alpha}$,
$(\exists x\in d\, A)^{(\exists\lambda\!\upharpoonright\!\alpha)}\equiv
(\exists x\in L_{\alpha}\, A)$.
\end{definition}

In what follows $F_{x,\lambda}$ denotes $F^{\Sigma_{1}}_{x,\lambda}$ when $\lambda\in R$, and $F^{\Sigma_{n}}_{x}$ when $\lambda=I$.

\begin{lemma}\label{lem:boundednessreg}{\rm (Boundedness)}
\\
Let $\lambda\in R^{+}$, $C\equiv (\exists x\in d\, A)$ and $C\not\in\{\exists x<\lambda \exists y<\lambda[\alpha<x \land P(\lambda,x,y)]: \alpha<\lambda\in R\}\cup\{\exists x<I[\alpha<x \land P_{I,n}(x)]:\alpha<I\}$.
Assume that $\mbox{{\rm rk}}(C)=\lambda=\mbox{{\rm rk}}_{L}(d)$.

\benu

\item\label{lem:boundednessregexi}

\[
(\mathcal{H},\Theta,\lambda,n)\vdash^{a}_{c}\Lambda, C \,\&\, a\leq b\in\mathcal{H}\cap\lambda
\Rightarrow (\mathcal{H},\Theta,\lambda,n)\vdash^{a}_{c}\Lambda,C^{(\exists\lambda\!\upharpoonright\! b)}
.\]
\item\label{lem:boundednessregfal}
\[
(\mathcal{H},\Theta,\kappa,n)\vdash^{a}_{c}\Lambda,\lnot C \,\&\, b\in\mathcal{H}\cap\lambda 
\Rightarrow (\mathcal{H},\Theta,\kappa,n)\vdash^{a}_{c}\Lambda,\lnot (C^{(\exists\lambda\!\upharpoonright\! b)})
.\]
\eenu

\end{lemma}

Though the following 
 Lemma \ref{th:Collapsingthmreg1}(Collapsing down to $I$)
is seen as in Lemma 5.22(Collapsing) of \cite{liftupZF},
we reproduce a proof of it since \cite{liftupZF} has not yet been published.

Recall that
\[
(\mathcal{H},\kappa,n)\vdash^{a}_{b}\Gamma
 :\Leftrightarrow 
(\mathcal{H},\emptyset,\kappa,n)\vdash^{a}_{b}\Gamma
\]

\begin{lemma}\label{th:Collapsingthmreg1}{\rm (Collapsing down to $I$)}\\

Suppose $\gamma\in\mathcal{H}_{\gamma,n}[\Theta]$ with $\Theta\subset\mathcal{H}_{\gamma,n}(\Psi_{I,n}\gamma)$, and 

\[
\Gamma\subset\Sigma^{\Sigma_{n+1}}(I)
\]

Then for $\hat{a}=\gamma+\omega^{I+a}$ 
\[
 (\mathcal{H}_{\gamma,n}[\Theta],I,n)\vdash^{a}_{I+1}\Gamma \Rightarrow
  (\mathcal{H}_{\hat{a}+1,n}[\Theta],I,n)\vdash^{\Psi_{I,n}\hat{a}}_{\Psi_{I,n}\hat{a}}\Gamma.
\]

\end{lemma}
\bprf

By induction on $a$. 
  
 First note that $\Psi_{I,n}\hat{a}\in\mathcal{H}_{\hat{a}+1,n}[\Theta]=\mathcal{H}_{\hat{a}+1,n}(\Theta)$
since 
$\hat{a}=\gamma+\omega^{I+a}\in\mathcal{H}_{\gamma,n}[\Theta]\subset\mathcal{H}_{\hat{a}+1,n}[\Theta]$
 by the assumption,
$\{\gamma,a\}\subset\mathcal{H}_{\gamma,n}[\Theta]$.
 
Assume $(\mathcal{H}_{\gamma,n}[\Theta][\Lambda],I,n)\vdash^{a_{0}}_{I+1}\Gamma_{0}$ with 
$\Lambda\subset\mathcal{H}_{\gamma,n}(\Psi_{I,n}\gamma)$.
Then by $\gamma\leq\hat{a}$, we have 
$\hat{a_{0}}\in \mathcal{H}_{\gamma,n}[\Theta][\Lambda]\subset \mathcal{H}_{\gamma,n}(\Psi_{I,n}\gamma)\subset\mathcal{H}_{\hat{a},n}(\Psi_{I,n}\hat{a})$.
This yields
that
\begin{equation}\label{eq:collapsethm}
a_{0}<a \Rightarrow \Psi_{I,n}\widehat{a_{0}}<\Psi_{I,n}\hat{a}
\end{equation}

Second observe that ${\sf k}^{E}_{\mathcal{K}}(\Gamma)\subset\mathcal{H}_{\gamma,n}[\Theta]\subset\mathcal{H}_{\hat{a}+1,n}[\Theta]$ by $\gamma\leq\hat{a}+1$.
 
 Third we have
 \begin{equation}\label{eq:Collapsingthm}
{\sf k}^{E}_{\mathcal{K}}(\Gamma)\subset\mathcal{H}_{\gamma,n}(\Psi_{I,n}\gamma)
\end{equation} 

\noindent
{\bf Case 1}.
First consider the case: $\Gamma\ni A\simeq\bigwedge\{A_{\iota} :\iota\in J\}$
\[
\infer[(\bigwedge)]{(\mathcal{H}_{\gamma,n}[\Theta],I,n)\vdash^{a}_{I+1}\Gamma}
{
 \{(\mathcal{H}_{\gamma,n}[\Theta\cup\{\mbox{{\rm rk}}_{L}(\iota)\}],I,n)\vdash^{a(\iota)}_{I+1}\Gamma, A_{\iota}:\iota\in J\}
 }\]
where $a(\iota)<a$ for any $\iota\in J$.

We claim that
\begin{equation}\label{eq:sigcollapsfal}
\forall\iota\in J(\mbox{{\rm rk}}_{L}(\iota)\in\mathcal{H}_{\gamma,n}(\Psi_{I,n}\gamma))
\end{equation}
Consider the case when $A\equiv\forall x\in b\, \lnot A^{\prime}$.
There are two cases to consider.
First consider the case when 
$J=\{d\}$ for the set
$d=\mu x\in b\, A^{\prime}$.
Then ${\sf k}^{E}_{\mathcal{K}}(A)={\sf k}(A)$, and
$\iota=d=(\mu x\in b\,  A^{\prime})\in\mbox{{\rm Hull}}_{\Sigma_{n}}^{I}({\sf k}(A))$, and 
$\mbox{{\rm rk}}_{L}(\iota)\in\mbox{{\rm Hull}}_{\Sigma_{n}}^{I}({\sf k}(A))\subset \mathcal{H}_{\gamma,n}(\Psi_{I,n}\gamma)$
by (\ref{eq:Collapsingthm}).
Otherwise 
we have $J=b$ and either $A\in\Sigma^{1}_{0}(\mathcal{K}^{+})$ and $b\in L_{\mathcal{K}^{+}}\cup\{L_{\mathcal{K}^{+}}\}$, or
$\mbox{{\rm rk}}_{L}(b)<I$.
In the second case we have $b\in{\sf k}(A)={\sf k}^{E}_{\mathcal{K}}(A)\subset\mathcal{H}_{\gamma,n}[\Theta]$.
In the first case each $\iota\in b$ has $L$-rank $\mbox{{\rm rk}}_{L}(\iota)<\mathcal{K}^{+}$.
On the other hand we have $\mathcal{K}^{+}\in\mathcal{H}_{\gamma,n}(\Psi_{I,n}\gamma)\cap I\subset\Psi_{I,n}\gamma$
 by $I>\mathcal{K}^{+}$.
Thus $\mbox{{\rm rk}}_{L}(\iota)<\Psi_{I,n}\gamma$.
In the second case we have $\mbox{{\rm rk}}_{L}(\iota)\leq\mbox{{\rm rk}}_{L}(b)\in \mathcal{H}_{\gamma,n}(\Psi_{I,n}\gamma)\cap I\subset\Psi_{I,n}\gamma$
by $\mbox{{\rm rk}}_{L}(b)<I$.

Hence (\ref{eq:sigcollapsfal}) was shown.

SIH  yields
\[
\infer[(\bigwedge)]
{(\mathcal{H}_{\hat{a}+1,n}[\Theta],I,n)\vdash^{\Psi_{I,n}\hat{a}}_{\Psi_{I,n}\hat{a}}\Gamma}
{
\{(\mathcal{H}_{\widehat{a(\iota)}+1,n}[\Theta\cup\{\mbox{{\rm rk}}_{L}(\iota)\}],I,n)\vdash^{\Psi_{I,n}\widehat{a(\iota)}}_{\Psi_{I,n}\widehat{a(\iota)}}
\Gamma, A_{\iota}:\iota\in J\}
}
\]
for $\widehat{a(\iota)}=\gamma+\omega^{I+a(\iota)}$, since $\Psi_{I,n}\widehat{a(\iota)}<\Psi_{I,n}\hat{a}$ by 
(\ref{eq:collapsethm}).
\\

\noindent
{\bf Case 2}.
Next consider the case for an $A\simeq\bigvee\{A_{\iota}:\iota\in J\}\in\Gamma$ and an $\iota\in J$
with $a(\iota)<a$ and $\mbox{{\rm rk}}_{L}(\iota)<I \Rightarrow \mbox{{\rm rk}}_{L}(\iota)<a$
\[
\infer[(\bigvee)]{(\mathcal{H}_{\gamma,n}[\Theta],I,n)\vdash^{a}_{I+1}\Gamma}
{(\mathcal{H}_{\gamma,n}[\Theta],I,n)\vdash^{a(\iota)}_{I+1}\Gamma,A_{\iota}}
\]
Assume $\mbox{{\rm rk}}_{L}(\iota)<I$. We show $\mbox{{\rm rk}}_{L}(\iota)<\Psi_{I,n}\hat{a}$.
By $\Psi_{I,n}\gamma\leq\Psi_{I,n}\hat{a}$,
 it suffices to show $\mbox{{\rm rk}}_{L}(\iota)<\Psi_{I,n}\gamma$.
 
Consider the case when $A\equiv\exists x\in b\, A^{\prime}$.
There are two cases to consider.
First consider the case when 
$J=\{d\}$ for the set
$d=\mu x\in b\, A^{\prime}$.
Then ${\sf k}^{E}_{\mathcal{K}}(A)={\sf k}(A)$, and
$\iota=d=(\mu x\in b\,  A^{\prime})\in\mbox{{\rm Hull}}_{\Sigma_{n}}^{I}({\sf k}(A))$, and 
$\mbox{{\rm rk}}_{L}(\iota)\in\mbox{{\rm Hull}}_{\Sigma_{n}}^{I}({\sf k}(A))\subset \mathcal{H}_{\gamma,n}(\Psi_{I,n}\gamma)$
by (\ref{eq:Collapsingthm}).
If $\mbox{{\rm rk}}_{L}(\iota)<I$, then $\mbox{{\rm rk}}_{L}(\iota)\in \mathcal{H}_{\gamma,n}(\Psi_{I,n}\gamma)\cap I\subset\Psi_{I,n}\gamma$.

Otherwise
we have $J=b$, and
either $A\in\Sigma^{1}_{0}(\mathcal{K}^{+})$ and $b\in L_{\mathcal{K}^{+}}\cup\{L_{\mathcal{K}^{+}}\}$, or
$b\in{\sf k}(A)={\sf k}^{E}_{\mathcal{K}}(A)\subset\mathcal{H}_{\gamma,n}[\Theta]$.
In the second case
we can assume that $\iota\in{\sf k}(A_{\iota})={\sf k}^{E}_{\mathcal{K}}(A_{\iota})\subset\mathcal{H}_{\gamma,n}[\Theta]$.
Otherwise set $\iota=0$.

In the first case each $\iota\in b$ has $L$-rank $\mbox{{\rm rk}}_{L}(\iota)<\mathcal{K}^{+}$.
On the other hand we have $\mathcal{K}^{+}\in\mathcal{H}_{\gamma,n}(\Psi_{I,n}\gamma)\cap I\subset\Psi_{I,n}\gamma$
 by $I>\mathcal{K}^{+}$.
Thus $\mbox{{\rm rk}}_{L}(\iota)<\Psi_{I,n}\gamma$.
In the second case we have 
$\mbox{{\rm rk}}_{L}(\iota)<\mbox{{\rm rk}}_{L}(b)\leq I$, and $\mbox{{\rm rk}}_{L}(\iota)\in \mathcal{H}_{\gamma,n}(\Psi_{I,n}\gamma)\cap I\subset\Psi_{I,n}\gamma$.

 SIH yields for $\widehat{a(\iota)}=\gamma+\omega^{I+a(\iota)}$
\[
\infer[(\bigvee)]
{(\mathcal{H}_{\hat{a}+1,n}[\Theta],I,n)\vdash^{\Psi_{I,n}\hat{a}}_{\Psi_{I,n}\hat{a}}}
{
(\mathcal{H}_{\widehat{a(\iota)}+1,n}[\Theta],I,n)\vdash^{\Psi_{I,n}\widehat{a(\iota)}}_{\Psi_{I,n}\widehat{a(\iota)}}\Gamma, A_{\iota}
}
\]
{\bf Case 3}.
Third consider the case for an $a_{0}<a$ and a $C$ with $\mbox{rk}(C)<I+1$.
\[
\infer[(cut)]{(\mathcal{H}_{\gamma,n}[\Theta],I,n)\vdash^{a}_{I+1}\Gamma}
{
(\mathcal{H}_{\gamma,n}[\Theta],I,n)\vdash^{a_{0}}_{I+1}\Gamma,\lnot C
&
(\mathcal{H}_{\gamma,n}[\Theta],I,n)\vdash^{a_{0}}_{I+1}C,\Gamma
}
\]
{\bf Case 3.1}. $\mbox{{\rm rk}}(C)<I$.

We have by (\ref{eq:Collapsingthm}) ${\sf k}^{E}_{\mathcal{K}}(C)\subset\mathcal{H}_{\gamma,n}(\Psi_{I,n}\gamma)$.
Proposition \ref{lem:rank}.\ref{lem:rank1} yields 
$\mbox{{\rm rk}}(C)\in\mathcal{H}_{\gamma,n}(\Psi_{I,n}\gamma)\cap I\subset\Psi_{I,n}\gamma\leq\Psi_{I,n}\hat{a}$.
By Proposition \ref{lem:rank}.\ref{prp:rksig3}
we see that $\{\lnot C,C\}\subset\Sigma^{\Sigma_{n+1}}(I)$.

SIH yields for $\widehat{a_{0}}=\gamma+\omega^{I+a_{0}}$
\[
\infer[(cut)]{(\mathcal{H}_{\hat{a}+1,n}[\Theta],I,n)\vdash^{\Psi_{I,n}\hat{a}}_{\Psi_{I,n}\hat{a}}\Gamma}
{
(\mathcal{H}_{\widehat{a_{0}}+1,n}[\Theta],I,n)\vdash^{\Psi_{I,n}\widehat{a_{0}}}_{\Psi_{I,n}\widehat{a_{0}}}\Gamma,\lnot C
&
(\mathcal{H}_{\widehat{a_{0}}+1,n}[\Theta],I,n)\vdash^{\Psi_{I,n}\widehat{a_{0}}
}_{\Psi_{I,n}\widehat{a_{0}}}C,\Gamma
}
\]
{\bf Case 3.2}.
$\mbox{{\rm rk}}(C)=I$.

Then $C\in\Sigma^{\Sigma_{n+1}}(I)$.
$C$ is either a sentence $\exists x<I[\alpha<x \land P_{I,n}(x)]$, 
or a sentence $\exists x\in L_{I}\, A(x)$ with
 ${\sf qk}(A)<I$.

In the first case we have 
$(\mathcal{H}_{\gamma+1,n}[\Theta],I,n)\vdash^{a_{0}}_{I+1}\Gamma$ by Reduction \ref{lem:reduction}.\ref{lem:reduction3},
 and IH yields the lemma.

Consider the second case.
From the right uppersequent, 
SIH yields for $\widehat{a_{0}}=\gamma+\omega^{I+a_{0}}$ and $\beta_{0}=\Psi_{I,n}\widehat{a_{0}}\in\mathcal{H}_{\widehat{a_{0}}+1,n}[\Theta]$
\[
(\mathcal{H}_{\widehat{a_{0}}+1,n}[\Theta],I,n)\vdash^{\beta_{0}}_{\beta_{0}} C,\Gamma
\]
Then by Boundedness \ref{lem:boundednessreg}.\ref{lem:boundednessregexi} and $\beta_{0}\in\mathcal{H}_{\widehat{a_{0}}+1,n}[\Theta]$,
we have
\[
(\mathcal{H}_{\widehat{a_{0}}+1,n}[\Theta],I,n)\vdash^{\beta_{0}}_{\beta_{0}}C^{(\exists I\!\upharpoonright\! \beta_{0})},\Gamma
\]
On the other hand we have by Boundedness \ref{lem:boundednessreg}.\ref{lem:boundednessregfal}
from the left uppersequent
\[
(\mathcal{H}_{\widehat{a_{0}}+1,n}[\Theta],I,n)\vdash^{a_{0}}_{\mu}\Gamma,\lnot (C^{(\exists I\!\upharpoonright\! \beta_{0})})
\]
Moreover we have
$\lnot (C^{(\exists I\!\upharpoonright\! \beta_{0})})\in\Sigma^{\Sigma_{n+1}}(I)$.
SIH yields for 
$\widehat{a_{0}}<\widehat{a_{1}}=\widehat{a_{0}}+1+\omega^{I+a_{0}}=\gamma+\omega^{I+a_{0}}+1+\omega^{I+a_{0}}<\gamma+\omega^{I+a}=\hat{a}$
and $\beta_{1}=\Psi_{I,n}\widehat{a_{1}}$
\[
(\mathcal{H}_{\widehat{a_{1}}+1,n}[\Theta],I,n)\vdash^{\beta_{1}}_{\beta_{1}}
\Gamma,\lnot C^{(\exists I\!\upharpoonright\! \beta_{0})}
\]
Now we have $\widehat{a_{i}}\in \mathcal{H}_{\widehat{a_{i}},n}(\Psi_{I,n}\hat{a})$ and $\widehat{a_{i}}<\hat{a}$ for $i<2$, and hence
$\beta_{0}=\Psi_{I,n}\widehat{a_{0}}<\beta_{1}=\Psi_{I,n}\widehat{a_{1}}<\Psi_{I,n}\hat{a}$. 
Therefore $\mbox{{\rm rk}}(C^{(\exists I\!\upharpoonright\! \beta_{0})})<\beta_{1}<\Psi_{I,n}\hat{a}$.

Consequently
\[
\infer[(cut)]{(\mathcal{H}_{\widehat{a_{1}}+1,n}[\Theta],I,n)\vdash^{\beta_{1}+1}_{\beta_{1}}\Gamma}
{
(\mathcal{H}_{\widehat{a_{1}}+1,n}[\Theta],I,n)\vdash^{\beta_{1}}_{\beta_{1}}\Gamma,
\lnot C^{(\exists I\!\upharpoonright\! \beta_{0})}
&
(\mathcal{H}_{\widehat{a_{0}}+1,n}[\Theta],I,n)\vdash^{\beta_{0}}_{\beta_{0}}C^{(\exists I\!\upharpoonright\! \beta_{0})},\Gamma
}
\]
Hence
$(\mathcal{H}_{\hat{a}+1,n},I,n)\vdash^{\Psi_{I,n}\hat{a}}_{\Psi_{I,n}\hat{a}}\Gamma$.
\\

\noindent
{\bf Case 4}. Fourth consider the case for an $a_{0}<a$
\[
\infer[(\mbox{{\bf F}})]{(\mathcal{H}_{\gamma,n}[\Theta],I,n)\vdash^{a}_{I+1}\Gamma}
{
(\mathcal{H}_{\gamma,n}[\Theta],I,n)\vdash^{a_{0}}_{I+1}\Lambda,\Gamma_{0}
}
\]
where $\Gamma=\Lambda\cup F"\Gamma_{0}$ and either  
$F=F^{\Sigma_{1}}_{x\cup\{\rho\}}$, $\Gamma_{0}\subset\Sigma_{1}$ for some $x$ and $\rho$,
or 
$F=F^{\Sigma_{n}}_{x}$, $\Gamma_{0}\subset\Sigma_{n}$ for an $x$.
Then $\Lambda\cup\Gamma_{0}\subset\Sigma_{n}$.
SIH yields the lemma.
\eprf

\subsection{Elimination of $\Pi^{1}_{1}$-indescribability}

In the subsection we eliminate inferences $(\mbox{{\bf Ref}}_{\mathcal{K}})$
for $\Pi^{1}_{1}$-indescribability.

For second-order sentences $\varphi$ on $L_{\pi}$ with parameters $A\subset L_{\pi}$ and ordinals $\alpha<\pi$, 
$\varphi^{(\alpha,\pi)}$ denotes the result of replacing second-order quantifiers $\exists X\subset L_{\pi},\forall X\subset L_{\pi}$ by
$\exists X\subset L_{\alpha},\forall X\subset L_{\alpha}$, resp., 
first-order quantifiers $\exists x\in L_{\pi},\forall x\in L_{\pi}$ by
$\exists x\in L_{\alpha},\forall x\in L_{\alpha}$, resp. and the parameters $A$ by $A\cap L_{\alpha}$.
For sequents $\Gamma$, $\Gamma^{(\alpha,\pi)}:=\{\varphi^{(\alpha,\pi)}:\varphi\in\Gamma\}$.

\begin{proposition}\label{prp:Mahloness}
Let $\Gamma\subset\Pi^{1}_{1}(\pi)$ for $\pi\in Mh_{n}^{\alpha}[\Theta]$.
Assume
\[
\exists\xi\in\mathcal{H}_{\xi,n}[\Theta\cup\{\pi\}](\pi)\cap\alpha
\forall\rho\in Mh_{n}^{\xi}[\Theta\cup\{\pi\}]
\bigvee(\Gamma^{(\rho,\pi)})
.\]
Then
$\bigvee(\Gamma)$ is true.
\end{proposition}
\bprf

By $\pi\in Mh_{n}^{\alpha}[\Theta]$ we have
$\pi\in M(Mh_{n}^{\xi}[\Theta\cup\{\pi\}])$ for any $\xi\in\mathcal{H}_{\xi,n}[\Theta\cup\{\pi\}](\pi)\cap\alpha$, cf. (\ref{eq:dfMh}).

Suppose the $\Sigma^{1}_{1}(\pi)$-sentence $\varphi:=\bigwedge(\lnot\Gamma):=\bigwedge\{\lnot\theta:\theta\in\Gamma\}$ is true.
Then 
the set $\{\rho<\pi: 
\varphi^{(\rho,\pi)}\}$ is club in $\pi$.

Hence for any $\xi\in\mathcal{H}_{\xi,n}[\Theta\cup\{\pi\}](\pi)\cap\alpha$ 
we can pick a $\rho\in Mh_{n}^{\xi}[\Theta\cup\{\pi\}]$ such that 
$ \varphi^{(\rho,\pi)}$.
\eprf

\[
\mathcal{H}_{\gamma,n}[\Theta]\vdash^{a}_{b}\Gamma
 :\Leftrightarrow 
(\mathcal{H}_{\gamma,n},\Theta,I,n)\vdash^{a}_{b}\Gamma.
\]

\begin{lemma}\label{lem:CollapsingthmKR}{\rm (Collapsing down to $\mathcal{K}$)}\\

Let $\gamma$ be
an ordinal such that
$
\gamma\in\mathcal{H}_{\gamma,n}
$.

Suppose for a finite set $\Theta$ of ordinals and an ordinal $a$
\[
\mathcal{H}_{\gamma,n}[\Theta]\vdash^{a}_{0}\Gamma
\]
where 
$\Gamma$ consists of sentences 
$\lnot\tau(B,\mathcal{K})$, $(B\cap C\neq\emptyset)$, $\forall\rho<\mathcal{K}\,\tau(B,\rho)$ for
 a $B\subset\mathcal{K}$ 
 with $B\in\mbox{{\rm Hull}}^{I}_{\Sigma_{1}}(\{\mathcal{K},\mathcal{K}^{+}\})$ and
 sets $C\in L_{\mathcal{K}+1}$ such that $C$ is a club subset of $\mathcal{K}$,
 and their subformulas:
\begin{equation}
\renewcommand{\theequation}{\ref{eq:thin}}
\tau(B,\rho):\Leftrightarrow \exists C\subset\rho
[(C \mbox{{ \rm is club}})^{\rho} \land
(B\cap C=\emptyset)]
\end{equation}
\addtocounter{equation}{-1}

Then for 
$\xi=\gamma+a$ 
\[
\forall\pi\in Mh_{n}^{\xi}[\Theta]\{
\models\Gamma^{(\pi,\mathcal{K})}
\}
.\]
which means that
$\bigvee (\Gamma^{(\pi,\mathcal{K})})$ is true for any $\pi\in Mh_{n}^{\xi}[\Theta]$.
\end{lemma}
\bprf

By induction on $a$. 
Let $\pi\in Mh_{n}^{\xi}[\Theta]$ and $\xi=\gamma+a$.
\\

\noindent
{\bf Case 1}.
First consider the case when the last inference is a $(\mbox{{\bf Ref}}_{\mathcal{K}})$:
we have $\{a_{\ell},a_{r}\}\subset\mathcal{H}_{\gamma,n}[\Theta]\cap a$
and $B\subset\mathcal{K}$ 
with $B\in\mbox{Hull}^{I}_{\Sigma_{1}}(\{\mathcal{K},\mathcal{K}^{+}\})$.

\[
\infer[(\mbox{{\bf Ref}}_{\mathcal{K}})]{\mathcal{H}_{\gamma,n}[\Theta]\vdash^{a}_{0}\Gamma}
{
\mathcal{H}_{\gamma,n}[\Theta]\vdash^{a_{\ell}}_{0}\Gamma,\lnot\tau(B,\mathcal{K})
&
  \mathcal{H}_{\gamma,n}[\Theta]\vdash^{a_{r}}_{0}\Gamma,\forall\rho<\mathcal{K}\, \tau(B,\rho)
}
\]

We have
$\xi_{r}:=\gamma+a_{r}\in\mathcal{H}_{\xi_{r},n}[\Theta](\pi)\cap\xi$ by $\xi_{r}\geq\gamma$ and $a_{r}<a$.
By Proposition \ref{prp:clshull}.\ref{prp:Mh3} with $\xi_{r}\in\mathcal{H}_{\xi_{r},n}[\Theta](\pi)$ we have
$\pi\in Mh_{n}^{\xi_{r}}[\Theta]$.
IH yields $\bigvee(\Gamma^{(\pi,\mathcal{K})})\lor \forall\rho<\pi\, \tau(B,\rho)$.

On the other hand we have $\xi_{\ell}:=\gamma+a_{\ell}\in\mathcal{H}_{\xi_{\ell},n}[\Theta](\pi)\cap\xi$.
By IH we have for any $\rho\in Mh_{n}^{\xi_{\ell}}[\Theta\cup\{\pi\}]\cap\pi$,
$\bigvee(\Gamma^{(\rho,\mathcal{K})})\lor \lnot\tau(B,\rho)$.
Hence we have
$\forall\rho\in Mh_{n}^{\xi_{\ell}}[\Theta\cup\{\pi\}]\cap\pi\{\bigvee(\Gamma^{(\rho,\mathcal{K})})\lor\bigvee(\Gamma^{(\pi,\mathcal{K})})\}$.
Proposition \ref{prp:Mahloness} yields $\bigvee(\Gamma^{(\pi,\mathcal{K})})$.
\\

\noindent
{\bf Case 2}.
Second consider the case when the last inference introduces a $\Pi^{1}_{1}(\mathcal{K})$-sentence $\lnot\tau(B,\mathcal{K})$
with  a $B\subset\mathcal{K}$ 
such that $B\in\mbox{Hull}^{I}_{\Sigma_{1}}(\{\mathcal{K},\mathcal{K}^{+}\})$.

\[
\infer[(\bigwedge)]{\mathcal{H}_{\gamma,n}[\Theta]\vdash^{a}_{0}\Gamma,\lnot\tau(B,\mathcal{K})}
{
\{
\mathcal{H}_{\gamma,n}[\Theta\cup\{\mbox{{\rm rk}}_{L}(C)\}]\vdash^{a(C)}_{0}\Gamma,
(C\not\subset\mathcal{K})\lor
\lnot (C\mbox{{ \rm is club}})^{\mathcal{K}}\lor
(B\cap C\neq\emptyset)
:
C\in L_{\mathcal{K}^{+}}
\}
}
\]
where $\forall C\in L_{\mathcal{K}^{+}}(a(C)\in\mathcal{H}_{\gamma,n}[\Theta\cup\{\mbox{{\rm rk}}_{L}(C)\}]\cap a)$ and 
$\lnot\tau(B,\mathcal{K})\simeq\bigwedge\{
(C\not\subset\mathcal{K})\lor
\lnot (C\mbox{{ \rm is club}})^{\mathcal{K}}\lor
(B\cap C\neq\emptyset)
:
C\in L_{\mathcal{K}^{+}}\}$.
For each $C$, $(C\not\subset\mathcal{K})\lor
\lnot (C\mbox{{ \rm is club}})^{\mathcal{K}}\lor
(B\cap C\neq\emptyset)$ is stratified with respec to $C$.

Let
\[
C_{\pi}:=\mu C\in L_{\pi^{+}}[(C\subset\pi)\land (C \mbox{{ \rm is club}})^{\pi} \land (B\cap C=\emptyset)]
\]
Then
$\lnot[(C_{\pi}\subset\pi)\land (C_{\pi} \mbox{{ \rm is club}})^{\pi} \land (B\cap C_{\pi}=\emptyset)]
\Rightarrow \lnot\tau(B,\pi)\equiv(\lnot\tau(B,\mathcal{K}))^{(\pi,\mathcal{K})}$.

We can assume that $(C_{\pi}\subset\pi)\land (C_{\pi} \mbox{{ \rm is club}})^{\pi}$.
Otherwise $(\lnot\tau(B,\mathcal{K}))^{(\pi,\mathcal{K})}$ and hence $\bigvee(\Gamma^{(\pi,\mathcal{K})})\lor(\lnot\tau(B,\mathcal{K}))^{(\pi,\mathcal{K})}$.

Let
\[
C=\{\gamma\in\mathcal{K}: \exists x,y<\mathcal{K}(\gamma=\pi \cdot x+y \land y\in C_{\pi}\cup\{0\})\}
\]
Then $C$ is an $L_{\mathcal{K}}$-definable club subset of $\mathcal{K}$, $C\in L_{\mathcal{K}+1}$, and
\\
$C\in J\cap\mbox{Hull}_{\Sigma_{1}}^{I}(\{\pi,\pi^{+},\mathcal{K},B\})\subset
\mbox{Hull}_{\Sigma_{1}}^{I}(\{\pi,\pi^{+},\mathcal{K},\mathcal{K}^{+}\})\subset
\mathcal{H}_{\gamma,n}[\Theta\cup\{\pi\}]$.
Hence 
$\mbox{{\rm rk}}_{L}(C)\in\mathcal{H}_{\gamma,n}[\Theta\cup\{\pi\}]$ and $a(C)\in\mathcal{H}_{\gamma,n}[\Theta\cup\{\pi\}]$.
By inversion
\[
\mathcal{H}_{\gamma,n}[\Theta\cup\{\pi\}]\vdash^{a(C)}_{0}\Gamma , C\not\subset\mathcal{K},
\lnot (C\mbox{{ \rm is club}})^{\mathcal{K}}, B\cap C\neq\emptyset
\]
Eliminate false sentences $C\not\subset\mathcal{K}$ and $\lnot (C\mbox{{ \rm is club}})^{\mathcal{K}}$
by Lemma \ref{lem:falsedelta0elim}.
\[
\mathcal{H}_{\gamma,n}[\Theta\cup\{\pi\}]\vdash^{a(C)}_{0}\Gamma , B\cap C\neq\emptyset
\]
IH yields for $\xi(C)=\gamma+a(C)$,
$\forall\rho\in Mh_{n}^{\xi(C)}[\Theta\cup\{\pi\}]\cap\pi\{\bigvee(\Gamma^{(\rho,\mathcal{K})})\lor(B\cap C\neq\emptyset)^{(\rho,\mathcal{K})}\}$,
where
$(B\cap C\neq\emptyset)^{(\rho,\mathcal{K})}\equiv(B\cap C_{\pi}\cap\rho\neq\emptyset)
\equiv((B\cap C\neq\emptyset)^{(\pi,\mathcal{K})})^{(\rho,\pi)}$.
Proposition \ref{prp:Mahloness} 
with $\xi(C)\in\mathcal{H}_{\xi(C),n}[\Theta\cup\{\pi\}](\pi)\cap\xi$
yields 
$\bigvee(\Gamma^{(\pi,\mathcal{K})})\lor(B\cap C\neq\emptyset)^{(\pi,\mathcal{K})}$,
and hence $\bigvee(\Gamma^{(\pi,\mathcal{K})})\lor(\lnot\tau(B,\mathcal{K}))^{(\pi,\mathcal{K})}$.
\\

\noindent
{\bf Case 3}. Third  consider the case : 
$\Gamma\ni (B\cap C\neq\emptyset)$ with 
$B\subset\mathcal{K}$, $B\in\mbox{Hull}^{I}_{\Sigma_{1}}(\{\mathcal{K},\mathcal{K}^{+}\})$ and
a club subset $C$ of $\mathcal{K}$.

\[
\infer[(\bigvee)]{\mathcal{H}_{\gamma,n}[\Theta]\vdash^{a}_{0}\Gamma}
{
\mathcal{H}_{\gamma,n}[\Theta]\vdash^{a_{0}}_{0}\Gamma, (d\in B)\land(d\in C)
 }
 \]
where $a_{0}<a$ and $d\in\mathcal{K}$.

Then $(B\cap C\neq\emptyset)^{(\pi,\mathcal{K})}\leftrightarrow (B\cap C\cap\pi\neq\emptyset)$ and $((d\in B)\land(d\in C))^{(\pi,\mathcal{K})}\leftrightarrow
(d\in(B\cap\pi))\land(d\in(C\cap\pi))$.
IH with Proposition \ref{prp:clshull}.\ref{prp:Mh3} yields the lemma.
\\

\noindent
{\bf Case 4}. Fourth consider the case : 
$\Gamma\ni ((d\in B)\land(d\in C))$ with 
$B\subset\mathcal{K}$, $B\in\mbox{Hull}^{I}_{\Sigma_{1}}(\{\mathcal{K},\mathcal{K}^{+}\})$ and
a club subset $C$ of $\mathcal{K}$.

\[
\infer[(\bigwedge)]{\mathcal{H}_{\gamma,n}[\Theta]\vdash^{a}_{0}\Gamma}
{
\mathcal{H}_{\gamma,n}[\Theta]\vdash^{a_{0}}_{0}\Gamma, d\in B
&
\mathcal{H}_{\gamma,n}[\Theta]\vdash^{a_{1}}_{0}\Gamma, d\in C
 }
 \]
where $a_{0},a_{1}<a$.

IH with Proposition \ref{prp:clshull}.\ref{prp:Mh3} yields the lemma.
\\

\noindent
{\bf Case 5}.
Fifth consider the case:
for a true literal $M\equiv(d\in B)$, $M\in\Gamma$,
where $B\subset\mathcal{K}$ such that either $B\in\mbox{Hull}^{I}_{\Sigma_{1}}(\{\mathcal{K},\mathcal{K}^{+}\})$, or
$B$ is a club subset of $\mathcal{K}$,
and $d\in\mathcal{K}$.
\[
\infer[(\bigwedge)]{\mathcal{H}_{\gamma,n}[\Theta]\vdash^{a}_{0}\Gamma}
{}
\]
Then $M^{(\pi,\mathcal{K})}\equiv(d\in(B\cap\pi))\in\Gamma^{(\pi,\mathcal{K})}$.

It suffices to show $d=\mbox{{\rm rk}}_{L}(d)<\pi$.
We have by (\ref{eq:controlder})
$d\in{\sf k}^{E}(d\in B)\cap\mathcal{K}\subset\mathcal{H}_{\gamma,n}\cap\mathcal{K}\subset\pi$
by $\pi\in Mh_{n}^{\xi}[\Theta]$, i.e., by $\mathcal{H}_{\xi,n}(\pi)\cap\mathcal{K}\subset\pi$.
\\

\noindent
{\bf Case 6}.
Sixth consider the case when the last inference introduces a sentence $\forall\rho<\mathcal{K}\,\tau(B,\rho)$.
\[
\infer[(\bigwedge)]{\mathcal{H}_{\gamma,n}[\Theta]\vdash^{a}_{0}\Gamma,\forall\rho<\mathcal{K}\,\tau(B,\rho)}
{
\{
\mathcal{H}_{\gamma,n}[\{\rho\}][\Theta]\vdash^{a(\rho)}_{0}\Gamma,\tau(B,\rho):
\rho<\mathcal{K}\}
}
\]
We have for any $\rho<\pi$ and $\xi(\rho)=\gamma+a(\rho)$, $\xi(\rho)\in\mathcal{H}_{\xi(\rho),n}[\Theta](\pi)$.
Proposition \ref{prp:clshull}.\ref{prp:Mh3} yields $\pi\in Mh_{n}^{\xi(\rho)}[\Theta]$.
By IH we have
$\forall\rho<\pi\{\bigvee(\Gamma^{(\pi,\mathcal{K})})\lor \tau(B,\rho)\}$, and hence
$(\bigvee(\Gamma)\lor \forall\rho<\mathcal{K}\,\tau(B,\rho))^{(\pi,\mathcal{K})}$
with $(\forall\rho<\mathcal{K}\,\tau(B,\rho))^{(\pi,\mathcal{K})}\equiv\forall\rho<\pi\,\tau(B,\rho)$.
\\

\noindent
{\bf Case 7}.
Seventh consider the case when the last inference introduces a sentence 
$\forall x\in c\,\varphi(x)\in\Gamma$ for $c\in L_{\mathcal{K}}$ and ${\sf k}^{E}(\varphi(x))<\mathcal{K}\,\&\, {\sf k}(\varphi(x))<\mathcal{K}^{+}$.
\[
\infer[(\bigwedge)]{\mathcal{H}_{\gamma,n}[\Theta]\vdash^{a}_{0}\Gamma}
{
\{
\mathcal{H}_{\gamma,n}[\{\mbox{{\rm rk}}_{L}(b)\}][\Theta]\vdash^{a(\rho)}_{0}\Gamma,\varphi(b):
b\in c\}
}
\]
Then $\gamma=\mbox{{\rm rk}}_{L}(c)\in{\sf k}^{E}(\Gamma)\cap\mathcal{K}$ and hence $\gamma<\pi$ as in {\bf Case 5}.
As in {\bf Case 6} we have by IH $\forall b\in c(\bigvee(\Gamma^{(\pi,\mathcal{K})})\lor\varphi(b))$ where
$\varphi(b)\equiv(\varphi(b))^{(\pi,\mathcal{K})}$.
Hence $\bigvee(\Gamma^{(\pi,\mathcal{K})})$.
\\

\noindent
{\bf Case 8}.
Eighth consider the case when the last inference introduces a sentence 
$\exists x\in c\,\varphi(c)\in\Gamma$ for $c\in L_{\mathcal{K}}$, $b\in c$ and ${\sf k}^{E}(\varphi(x))<\mathcal{K}\,\&\, {\sf k}(\varphi(x))<\mathcal{K}^{+}$.
\[
\infer[(\bigvee)]{\mathcal{H}_{\gamma,n}[\Theta]\vdash^{a}_{0}\Gamma}
{
\mathcal{H}_{\gamma,n}[\Theta]\vdash^{a_{0}}_{0}\Gamma,\varphi(b)
}
\]
As in {\bf Case 7} we see $\mbox{{\rm rk}}_{L}(c)<\pi$.
IH with Proposition \ref{prp:clshull}.\ref{prp:Mh3} yields $\bigvee(\Gamma^{(\pi,\mathcal{K})})\lor\varphi(b)$, and 
$\bigvee(\Gamma^{(\pi,\mathcal{K})})$.
\\

\noindent
{\bf Case 9}.
Ninth consider the case when the last inference is an $(\mbox{{\bf F}})$
where
either $F=F^{\Sigma_{1}}_{x\cup\{\lambda\}}$ for a $\lambda\in R$
or 
$F=F^{\Sigma_{n}}_{x}$.

In each case if $A\in rng(F)$ for an $A\in\Gamma$,
then we claim
$F"A\equiv A$.
Suppose $x=F^{\Sigma_{1}}_{x\cup\{\mathcal{K}^{+}\}}(\mathcal{K}^{+})\leq \mbox{{\rm rk}}_{L}(B)<\mathcal{K}^{+}$
for the set $B\in\mbox{{\rm Hull}}^{I}_{\Sigma_{1}}(\{\mathcal{K},\mathcal{K}^{+}\})$.
However by $x>\mathcal{K}$ we have 
$\mbox{{\rm rk}}_{L}(B)\in\mbox{{\rm Hull}}^{I}_{\Sigma_{1}}(\{\mathcal{K},\mathcal{K}^{+}\})\cap\mathcal{K}^{+}\subset\mbox{{\rm Hull}}^{I}_{\Sigma_{1}}(x\cup\{\mathcal{K}^{+}\})\cap\mathcal{K}^{+}\subset x$.
Hence this is not the case.

IH yields the assertion.
\eprf

Collapsing down to $\mathcal{K}$ \ref{lem:CollapsingthmKR}
yields the following Theorem \ref{cor:CollapsingthmK}.

\begin{theorem}\label{cor:CollapsingthmK}{\rm (Elimination of $(\mbox{{\bf Ref}}_{\mathcal{K}})$)}\\

Let $\gamma\in\mathcal{H}_{\gamma,n}$, $B\subset\mathcal{K}$,
and $B\in\mbox{{\rm Hull}}^{I}_{\Sigma_{1}}(\{\mathcal{K},\mathcal{K}^{+}\})$. 

\[
[ \mathcal{H}_{\gamma,n}\vdash^{a}_{0}\lnot\tau(B,\mathcal{K})]
 \Rightarrow
[\lnot\tau(B,\pi) \mbox{ {\rm is true}}]
\]
for any $\pi\in Mh_{n}^{\xi}$ with $\xi=\gamma+a$.

\end{theorem}

\section{Proof of Theorem \ref{th:mainthK}}\label{sect:prfmainthm}

Let $\varphi$ be a $\Sigma^{1}_{2}$-sentence, and
assume that ${\sf ZF}$
proves the sentence
\[
\forall\mathcal{K}[(\mathcal{K} \mbox{ {\rm is a weakly compact cardinal}})  \to \varphi^{V_{\mathcal{K}}}]
.\]

Under $V=L$, $V_{\sigma}=L_{\sigma}$ for any inaccessible cardinals $\sigma$, and 
we have $\forall\mathcal{K}[(\mathcal{K} \mbox{ {\rm is a weakly compact cardinal}}) \to \varphi^{L_{\mathcal{K}}}]$.
Hence $\mbox{{\rm T}}(\mathcal{K},I)\vdash \varphi^{L_{\mathcal{K}}}$.
By Proposition \ref{prp:Jensen} we can assume that the sentence
(`$\mathcal{K} \mbox{ {\rm is uncountable regular}' } \to\varphi^{L_{\mathcal{K}}}$) is of the form `$\exists B\subset\mathcal{K}(S^{\varphi}(B)\cap\mathcal{K}$ is stationary in $\mathcal{K}$)'.

Let $B:=\mu B\subset\mathcal{K}(S^{\varphi}(B)\cap\mathcal{K}\mbox{ is stationary in }\mathcal{K})\in\mbox{{\rm Hull}}^{I}_{\Sigma_{1}}(\{\mathcal{K},\mathcal{K}^{+}\})$.

In what follows work in an intuitionistic fixed point theory $\mbox{FiX}^{i}(\mbox{{\sf ZFLK}}_{n})$
over
$\mbox{{\sf ZFLK}}_{n}={\sf ZF}+(V=L)+(\mathcal{K}\in Mh_{n}^{\omega_{n}(I+1)})$
for a sufficiently large $n<\omega$.
By Embedding \ref{th:embedreg} pick an $m<\omega$ such that 
$(\mathcal{H}_{0,n},I,n)\vdash^{I\cdot (m-1)}_{I+m-1}\lnot\tau(B,\mathcal{K})$.
By Predicative Cut-elimination \ref{lem:predcereg} we have
$(\mathcal{H}_{0,n},I,n)\vdash^{\omega_{m-2}(I\cdot (m-1))}_{I+1}\lnot\tau(B,\mathcal{K})$.

Then by Collapsing down to $I$ \ref{th:Collapsingthmreg1}
we have for $a=\omega_{m}(I+1)$ and 
$b=\Psi_{I,n}a$,
$(\mathcal{H}_{a,n},I,n)\vdash^{b}_{b}\lnot\tau(B,\mathcal{K})$.
Again by Predicative Cut-elimination \ref{lem:predcereg} 
we have 
$(\mathcal{H}_{a,n},I,n)\vdash^{\varphi bb}_{0}\lnot\tau(B,\mathcal{K})$.

Elimination of $(\mbox{{\bf Ref}}_{\mathcal{K}})$ \ref{cor:CollapsingthmK} 
yields
$\lnot\tau(B,\pi)$
for any $\pi\in Mh_{n}^{\xi}$ with
$\xi=a+\varphi bb\in\mathcal{H}_{\xi,n}(\mathcal{K})\cap\omega_{m+1}(I+1)$.

Proposition \ref{prp:Mahloness} with $\mathcal{K}\in Mh_{n}^{\omega_{m+1}(I+1)}$ yields $\lnot\tau(B,\mathcal{K})$,
 and hence $S^{\varphi}(B)\cap\mathcal{K}$ is stationary in $\mathcal{K}$.
 Since the whole proof is formalizable in $\mbox{FiX}^{i}(\mbox{{\sf ZFLK}}_{n})$, 
we conclude
$\mbox{FiX}^{i}(\mbox{{\sf ZFLK}}_{n})\vdash  \varphi^{V_{\mathcal{K}}}$.
Finally Theorem \ref{th:consvintfix} yields
$\mbox{{\sf ZFLK}}_{n}\vdash \varphi^{V_{\mathcal{K}}}$.
Therefore $\varphi^{V_{\mathcal{K}}}$ follows from $\theta_{n}(\mathcal{K}) :\Leftrightarrow \mathcal{K}\in Mh_{n}^{\omega_{n}(I+1)}$ over ${\sf ZF}+(V=L)$.
Thus Theorem \ref{th:mainthK}.\ref{th:mainthK2} was shown.

Since the least weakly inaccessible cardinal $I_{0}$ is below the least weakly Mahlo cardinal,
\[
{\sf ZF}+\mathbb{K}\vdash \varphi^{V_{I_{0}}}
\Rightarrow
{\sf ZF}+\{\exists\mathcal{K}\, \theta_{n}(\mathcal{K}) : n<\omega\}\vdash \varphi^{V_{I_{0}}}
\]
 for any first-order sentence $\varphi$, etc.

This completes a proof of Theorem \ref{th:mainthK}.

\end{document}